\documentclass[11pt]{amsart}

\usepackage{epigamath}

\usepackage[english]{babel}

\numberwithin{equation}{section}

\usepackage[shortlabels]{enumitem}

\usepackage[T1]{fontenc}
\usepackage{etex}
\usepackage{amsxtra}
\usepackage{amscd}
\usepackage{amsthm}
\usepackage{amsfonts}
\usepackage{tikz-cd}
\usepackage{mathrsfs}
\usepackage{subfiles}
\usepackage{morefloats}
\usepackage[utf8]{inputenc}
\usepackage{scalerel}
\usepackage{stackengine}
\usepackage{microtype}
\usepackage{cleveref}

\theoremstyle{plain}
\newtheorem{theorem}{Theorem}[section]
\newtheorem{proposition}[theorem]{Proposition}
\newtheorem{lemma}[theorem]{Lemma}
\newtheorem{corollary}[theorem]{Corollary}

\theoremstyle{definition}
\newtheorem{definition}[theorem]{Definition}
\newtheorem{example}[theorem]{Example}

\theoremstyle{remark}
\newtheorem{remark}[theorem]{Remark}

\AddToHook{env/proposition/begin}{\crefalias{theorem}{proposition}}
\AddToHook{env/lemma/begin}{\crefalias{theorem}{lemma}}
\AddToHook{env/corollary/begin}{\crefalias{theorem}{corollary}}
\AddToHook{env/definition/begin}{\crefalias{theorem}{definition}}
\AddToHook{env/example/begin}{\crefalias{theorem}{example}}
\AddToHook{env/remark/begin}{\crefalias{theorem}{remark}}

\newcommand\bc{\mathbb C}

\newcommand\bp{\mathbb P}

\newcommand\br{\mathbb R}

\newcommand\bv{\mathbb V}

\newcommand\bz{\mathbb Z}

\DeclareMathOperator\spec{Spec}

\newcommand*{\shom}{\mathscr{H}\kern -.5pt om}
\newcommand*{\stor}{\mathscr{T}\kern -.5pt or}
\newcommand*{\sext}{\mathscr{E}\kern -.5pt xt}

\DeclareMathOperator\id{id}

\renewcommand\hom{Hom}
\DeclareMathOperator\coker{coker}

\DeclareMathOperator\rank{rank}

\DeclareMathOperator\spn{span}
\DeclareMathOperator\im{im}
\DeclareMathOperator\End{End}

\renewcommand\sp{\mathrm{Sp}}

\DeclareMathOperator\tr{tr}

\DeclareMathOperator\res{Res}

\DeclareMathOperator\Mod{Mod}

\DeclareMathOperator\gl{GL}

\renewcommand\sl{\mathrm{SL}}
\newcommand\psl{\mathrm{PSL}}

\DeclareMathOperator\so{SO}
\DeclareMathOperator\su{SU}
\renewcommand\u{\mathrm{U}}

\DeclareMathOperator\ad{ad}

\DeclareMathOperator\gr{Gr}

\DeclareMathOperator\pardeg{par-deg}
\DeclareMathOperator\Rep{Rep}

\newcommand\lf{\mathrm{lf}}
\newcommand\tors{\mathrm{tors}}
\newcommand{\lra}{\longrightarrow}

\newcommand\sslash{\mathbin{/\mkern-3.5mu/}}

\setlist[enumerate,1]{label={\rm(\arabic*)}, ref={\rm\arabic*}}
\newcommand{\supth}[1]{\ensuremath{#1^{\mathrm{th}}}}

\EpigaVolumeYear{10}{2026} \EpigaArticleNr{23} \ReceivedOn{April 23, 2025}
\InFinalFormOn{May 7, 2026}
\AcceptedOn{July 1, 2026}

\title{Minimal energy local systems on curves}
\titlemark{Minimal energy local systems on curves}

\author{Charlie Wu}
\address{Department of Mathematics, University of Toronto, 40 St.\ George Street, Toronto, ON M5S 2E4, Canada}
\email{charliec.wu@mail.utoronto.ca}

\authormark{C.~Wu}

\AbstractInEnglish{Let $\Sigma_{g,n}$ be an orientable topological surface of genus $g$ with $n$ punctures. When $g = 0$, Deroin and Tholozan studied the class of \textit{supra-maximal} representations $\pi_1(\Sigma_{0,n})\to \psl_2(\br)$, and they showed that the supra-maximal representations form a compact component of a real relative character variety. We study a collection of rank $r$ local systems on $\Sigma_{g,n}$ that we call \textit{of minimal energy}. These are generalizations of supra-maximal representations, and underlie polarizable complex variations of Hodge structure for any choice of complex structure on $\Sigma_{g,n}$. Like the supra-maximal representations, the minimal energy local systems form compact components of relative character varieties of real forms of $\gl_r(\bc)$.

  We show that when the local monodromy data around the punctures is chosen to be unitary and generic, and the relative character variety is nonempty, these minimal energy local systems always exist. When $g > 0$, we show that the minimal energy local systems come from unitary representations of $\pi_1(\Sigma_{g,n})$. If $g = 0$, we show that they do not always come from unitary representations, and we study their structure in general.}


\begin{document}



\maketitle

\begin{prelims}

\DisplayAbstractInEnglish

\bigskip

\DisplayKeyWords

\medskip

\DisplayMSCclass

\end{prelims}


\newpage

\setcounter{tocdepth}{1}

\tableofcontents


\section{Introduction}\label{sec:intro}

Let $\Sigma_{g,n}$ be an orientable topological surface of genus $g$ and $n$ punctures, and let $G \subseteq \gl_r(\bc)$ be a reductive algebraic group over $\br$ or $\bc$. Let $C = (C_1, \dots, C_n)$ be an $n$-tuple of conjugacy classes of $G$. We let $\chi(\pi_1(\Sigma_{g,n}), G, C)$ be the \textit{relative character variety} of $\Sigma_{g,n}$. This is the space of conjugacy classes (\textit{i.e.}, isomorphism classes) of semisimple rank $r$ representations of the fundamental group of $\Sigma_{g,n}$ valued in $G$ such that the homotopy class of a simple closed loop $\gamma_i$ around some puncture $x_i$ lies in $C_i$.

More precisely, we let the \textit{relative representation variety} $\Rep(\pi_1(\Sigma_{g,n}), G, C)$ be given by
\begin{align*}
	\Rep\left(\pi_1\left(\Sigma_{g,n}\right), G, C\right) \subseteq  \prod_{i = 1}^{2g} G \times \prod_{i = 1}^n \overline{C}_i, 
\end{align*}
where $\overline{C}_i$ is the Zariski closure of $C_i$ in $G$, and $\Rep(\pi_1(\Sigma_{g,n}), G, C)$ is cut out by the relations 
\begin{align*}
	\prod_{i = 1}^{g} [A_{2i-1},A_{2i}] = \prod_{i = 1}^n B_i
\end{align*}
for $A_i \in G$ and $B_i \in \overline{C}_i$. Since the $\overline{C}_i$ are algebraic subvarieties of $G$, $\Rep(\pi_1(\Sigma_{g,n}), G, C)$ inherits the structure of an algebraic variety. Then, $G$ acts on $\Rep(\pi_1(\Sigma_{g,n}), G, C)$ by simultaneous conjugation of the $A_i$ and $B_i$, and we define $\chi(\pi_1(\Sigma_{g,n}), G, C)$ to be $\Rep(\pi_1(\Sigma_{g,n}), G, C)\sslash G$, where the quotient is a GIT quotient. This GIT quotient is also the topological quotient of $\Rep(\pi_1(\Sigma_{g,n}), G,C)$ by the conjugation action. Throughout, we work with conjugacy classes that are semisimple (so the $C_i$ are Zariski closed).

\subsection{Main results}
We study a collection of local systems (equivalently, $\pi_1(\Sigma_{g,n})$-representations) in $\chi(\pi_1(\Sigma_{g,n}), \gl_r(\bc), C)$ that we call \textit{of minimal energy}. These local systems are defined in \Cref{def: definition_minimal_energy}, and they come from polarizable complex variations of Hodge structure ($\bc$-VHS) satisfying a Hodge-theoretic vanishing condition. This definition \textit{a priori} depends on a choice of complex structure for $\Sigma_{g,n}$, but we show in \Cref{prop: top_dimension_C-VHS} that this choice turns out to be auxiliary and does not matter.

Therefore, we give $\Sigma_{g,n}$ a complex structure and work with smooth algebraic curves. Let $X$ be a smooth proper curve of genus $g$ over the complex numbers, and let $D = x_1 + \dots + x_n$ be a reduced effective divisor on $X$ (\textit{i.e.}, a collection of distinct points $\{x_1,\dots, x_n\})$ in $X$) so that topologically, $\Sigma_{g,n}$ is homeomorphic to $X \setminus D$. Suppose that the conjugacy classes $C_1,\dots, C_n \subseteq \gl_r(\bc)$ are \emph{elliptic}, meaning that they contain a unitary matrix. As before we let $C = (C_1,\dots, C_n)$. Since these conjugacy classes are semisimple (since they preserve a Hermitian metric), they are Zariski closed.

The following results will be proven in \Cref{sec:minimal energy,sec:genus zero,sec:consequences}.

We work with $\chi(\pi_1(X \setminus D), \gl_r(\bc),C)$ when the relative character variety is smooth. This condition is guaranteed provided that every representation $\rho \in \chi(\pi_1(X \setminus D), \gl_r(\bc),C)$ is irreducible, and we can ensure that every representation in $\chi(\pi_1(X \setminus D), \gl_r(\bc),C)$ is irreducible by picking our conjugacy classes generically. The precise condition that guarantees this is spelled out in \Cref{def: definition_generic_distinct};  it is a numerical condition on the eigenvalues of the conjugacy classes that can be easily checked by hand.

\begin{theorem}\label{thm: deform_minimal_energy}
If $\chi(\pi_1(X\setminus D), \gl_r(\bc),C)$ is nonempty and smooth, then there is a minimal energy local system in every connected component of $\chi(\pi_1(X \setminus D), \gl_r(\bc), C)$.
\end{theorem}

The methods used to prove \Cref{thm: deform_minimal_energy} allow us conclude in \Cref{prop: top_dimension_C-VHS} that the minimal energy local systems form almost all of the $\bc$-VHS lying in $\chi(\pi_1(X \setminus D), \gl_r(\bc),C)$. More precisely, we show that for each connected component of $\chi(\pi_1(X \setminus D),\gl_r(\bc),C)$, the locus of $\bc$-VHS in that connected component contains the locus of minimal energy local systems as the \emph{unique} connected component of the $\bc$-VHS locus of largest dimension in that connected component.

If $X$ has positive genus and our local monodromy data is generic, a minimal energy local system must come from a $\bc$-VHS with a one-step Hodge filtration. That is, these minimal energy local systems come from unitary representations. This is proven in \Cref{prop: minimal_energy_genus_0}. If the genus of $X$ is zero, then the minimal energy local systems do not have to come from a $\bc$-VHS with a one-step Hodge filtration. In \Cref{ex: two_graded_pieces_example}, for any rank $r$, we find an example of a minimal energy local system on $\bp^1 \setminus D$ whose associated $\bc$-VHS has a two-step Hodge filtration. Moreover, we show that when the number of punctures $n$ is very large compared to $r$, this is the maximum number of pieces that can appear in the Hodge filtration of a $\bc$-VHS associated to a minimal energy local system.

\begin{theorem}\label{thm: bounds_C-VHS}
Let $\mathbb{V}$ be an irreducible minimal energy local system of rank $r$ on $\bp^1 \setminus \{x_1,\dots, x_n\}$ with generic unitary local monodromy $($in the sense of \Cref{def: definition_generic_distinct}\,$)$. If
\begin{align*}
	n > 5r^2-25, 
\end{align*}
then $\mathbb{V}$ comes from a $\bc$-VHS with at most two steps in its Hodge filtration.
\end{theorem}
For $r \leq 2$, we have that $5r^2 - 25 \leq -5$, so the condition on the number of punctures $n$ is vacuous because the number of punctures is nonnegative. Indeed, every rank $2$ $\bc$-VHS can have only at most two steps in its Hodge filtration, so \Cref{thm: bounds_C-VHS} is uninteresting. However, when $r \geq 3$ (so that $n \geq 20$), \Cref{thm: bounds_C-VHS} places nontrivial conditions on the $\bc$-VHS. This bound is far from sharp, however. For example, a rank $3$ minimal energy local system on $\bp^1 \setminus \{x_1,\dots, x_n\}$ can have at most two steps in its Hodge filtration as soon as $n \geq 4$.

\subsection{History and motivation}

This paper and its motivations are heavily inspired by work of Benedetto--Goldman \cite{Benedetto_Goldman}, Deroin--Tholozan \cite{Deroin_Tholozan_super_maximal}, Mondello \cite{Mondello}, Tholozan--Toulisse \cite{Tholozan_Toulisse}, and Feng--Zhang \cite{Feng_Zhang}. These papers studied and produced interesting compact components of relative character varieties of punctured spheres for various real Lie groups $G$ and conjugacy classes $C = (C_1,\dots, C_n)$.

\begin{figure}[htbb]
  \centering
\includegraphics[width=0.5\linewidth]{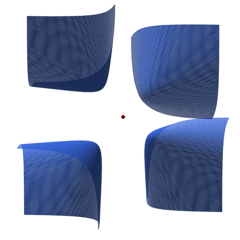}
  \caption{$\chi(\pi_1(\Sigma_{0,4}), \sl_2(\br), C)$ with $\tr(C_1) = \tr(C_2)= \tr(C_3) = 3/2$ and $\tr(C_4) = -3/2$. The tiny red component is compact, and the blue components are noncompact.}\label{fig:BG variety}
\end{figure}

Benedetto and Goldman first studied the case of a four-punctured sphere with $G = \sl_2(\br) \cong \su(1,1)$ and $G = \su(2)$. They showed that $\chi(\pi_1(\Sigma_{0,4}), \sl_2(\bc), C)$ admits a real form with a compact component whenever the traces of the $C_i$ lie in $(-2,2)$. This compact component  lies in either $\chi(\pi_1(\Sigma_{0,4}), \sl_2(\br), C)$ or $\chi(\pi_1(\Sigma_{0,4}), \su(2), C)$. The group $\su(2)$ is compact, so this is not so surprising. The group $\sl_2(\br)$ is not compact, but nevertheless Benedetto and Goldman surprisingly find examples of conjugacy classes $C = (C_1,\dots, C_4)$ where $\chi(\pi_1(\Sigma_{0,4}), \sl_2(\br), C)$ has a compact component. More precisely, they show that $\chi(\pi_1(\Sigma_{0,4}), \sl_2(\br), C)$ has a compact connected component when the conjugacy classes $(C_1,\dots, C_4)$ satisfy a certain explicit inequality in the traces of the $C_i$. When the inequality is not satisfied, the compact component shows up in the $\su(2)$-character variety. We point the reader to \cite[Theorem 3.12]{Cantat_Loray} for the explicit bound. Remarkably, Benedetto and Goldman found these components by graphing them on a computer (see \Cref{fig:BG variety}).

Deroin--Tholozan \cite[Corollary 3]{Deroin_Tholozan_super_maximal} then generalized the work of Benedetto--Goldman to spheres with arbitrarily many punctures. Using methods from geometric topology, they construct the class of \emph{supra-maximal} representations $\rho\colon \pi_1(\Sigma_{0,n}) \to \sl_2(\br)$ and show that the collection of supra-maximal representations forms a compact component of $\chi(\pi_1(\Sigma_{0,n}), \sl_2(\br), C)$ for particular choices of elliptic conjugacy classes $C = (C_1,\dots, C_n)$. Along the way, they show that these supra-maximal representations have very interesting algebraic and geometric properties. They show that the homotopy class of every simple closed curve in $\Sigma_{0,n}$ is sent to a nonhyperbolic element of $\sl_2(\br)$, see \cite[Theorem 2]{Deroin_Tholozan_super_maximal}, and for any supra-maximal representation $\rho$ and any complex structure on $\Sigma_{0,n}$, there is a unique $\rho$-equivariant holomorphic map from the universal cover of $\Sigma_{0,n}$ to the symmetric space of $\sl_2(\br)$ (the upper half-plane). Related work of Mondello \cite[Corollary 3.20]{Mondello} studies $\chi(\pi_1(\Sigma_{g,n}), \sl_2(\br), C)$  for general $g$, and Mondello counts the number of components of $\chi(\pi_1(\Sigma_{g,n}), \sl_2(\br), C)$ and recovers the compact components. Mondello also shows that the compact component is \emph{unique} when it exists.

Feng--Zhang and Tholozan--Toulisse then found analogous compact components of real relative character varieties. Tholozan--Toulisse found compact components of $\chi(\pi_1(\Sigma_{0,n}), G, C)$ for the Hermitian Lie groups $G = \su(p,q)$, $\sp(2, \br)$, and $\so(2n)$. Using nonabelian Hodge theory, see \cite[Theorem 1]{Tholozan_Toulisse}, they recover similar algebraic and geometric properties of representations in their compact component; see \cite[Theorem 2]{Tholozan_Toulisse}. Their methods explicitly construct the Higgs bundles corresponding to these representations in the compact components.

Feng--Zhang found a compact component for $\chi(\pi_1(\Sigma_{0,n}), \so_o(2,q), C)$, where $\so_o(2,q)$ is the identity component of $\so(2,q)$; see \cite[Theorem 1.1(a)]{Feng_Zhang}. Their work covers the last collection of Hermitian Lie groups left open from the work of Tholozan--Toulisse. The representations in this component discovered by Feng--Zhang have the same good algebraic and geometric properties found in the components of Deroin--Tholozan and Tholozan--Toulisse; see \cite[Theorem 1.1(b), (c)]{Feng_Zhang}. Like the work of Tholozan--Toulisse, their methods involve explicitly constructing the Higgs bundles corresponding to these representations in the compact components.

We give a new interpretation of these compact components. Like the work of Feng--Zhang and Tholozan--Toulisse, we use nonabelian Hodge theory to study the relative character variety $\chi(\pi_1(\Sigma_{g,n}), G, C)$. However, we find compact components of $\chi(\pi_1(\Sigma_{g,n}), G, C)$ for $G = \u(p,q)$ without writing down the Higgs bundles. Instead, we use deformation-theoretic techniques in the moduli space of parabolic Higgs bundles. One major upside to our methods is that we do not need to work with the conjugacy classes in a strong way. Under mild hypotheses on the conjugacy classes (to ensure that $\chi(\pi_1(\Sigma_{g,n}), \gl_r(\bc), C)$ is nonempty and smooth), there are compact components of $\chi(\pi_1(\Sigma_{g,n}), \u(p,q), C)$ or $\chi(\pi_1(\Sigma_{g,n}), \su(p,q), C)$ (whose points are the minimal energy local systems) for some pair of nonnegative integers $(p,q)$ such that $p + q = r$. In particular, we recover the compact components found by Tholozan--Toulisse and Feng--Zhang (see \Cref{rem:old examples}).

Along the way, we show that these compact components (which consist of minimal energy local systems) are made up of \emph{universal} $\bc$-VHS. That is, they underlie complex variations of Hodge structure for any choice of complex structure on $\Sigma_{g,n}$. Since the unitary group is compact and all unitary local systems on $\Sigma_{g,n}$ come from $\bc$-VHS, the results are not surprising when the minimal energy local systems correspond to unitary local systems. However, we show that there are often minimal energy local systems in $\chi(\pi_1(\Sigma_{g,n}), \su(p,q), C)$ when $\su(p,q)$ is noncompact (that is, when $p$ and $q$ are both nonzero). This gives us both interesting compact components of $\chi(\pi_1(\Sigma_{g,n}), \su(p,q), C)$ and a class of interesting universal $\bc$-VHS on $\Sigma_{g,n}$.

\subsection{Structure of the paper}
Our methods heavily rely on nonabelian Hodge theory, which gives us a correspondence between local systems on $X \setminus D$ and parabolic Higgs bundles on $X$ (with parabolic divisor $D$). In \Cref{sec:parabolic bundles,sec:VHS} we review the notions of parabolic bundles and $\bc$-VHS. In particular, we discuss how to take a $\bc$-VHS on $X \setminus D$ and produce a parabolic Higgs bundle. 

Experts may wish to skip to \Cref{sec:minimal energy}, where we introduce the notion of minimal energy and prove \Cref{thm: deform_minimal_energy} and \Cref{prop: minimal_energy_genus_0}. In \Cref{sec:genus zero} we prove \Cref{thm: bounds_C-VHS}, and in \Cref{sec:consequences} we give \Cref{ex: two_graded_pieces_example} and give some applications of minimal energy local systems to Gromov--Witten theory.

\subsection*{Acknowledgments}
I thank my advisor Daniel Litt. He suggested that I study these minimal energy local systems, and this project would not be possible without his encouragement and the many extremely helpful conversations I have had with him. I would also like to thank Simon Xu for the discussions I have had with him and for entertaining my many questions about Hodge theory. I am grateful to Samuel Bronstein and Arnaud Maret for their insightful questions and comments. I also thank Brian Collier as well as Bronstein and Maret for bringing the work of Tholozan and Toulisse to my attention. I would also like to thank the anonymous referees for their very helpful and detailed comments, and for pointing out some mistakes in previous versions of the paper.

\section{Parabolic bundles}\label{sec:parabolic bundles}

We review the notions of parabolic bundles and parabolic sheaves. We use the notation from \cite[Section 2]{Landesman-Litt-very_general_curves}. Parabolic bundles are, roughly speaking, bundles carrying extra structure that keeps track of local monodromy data.

\begin{definition}\label{def: definition_parabolic_bundle}
Let $E$ be a vector bundle on a smooth curve $X$, and let $D = x_1 + \dots + x_n$ be a reduced effective divisor on $X$ (that is, $D$ is a collection of points). A \emph{parabolic bundle} $(E, \{E_j^i\} ,\{\alpha_j^i\})$ on $(X,D)$ is the data
\begin{enumerate}
	\item a strictly decreasing filtration $E_{x_j} = E_j^1 \supset E_j^2 \supset \dots \supset E_j^{n_j+1} = 0$,
	\item a sequence of real numbers $0 \leq \alpha_j^1 < \alpha_j^2 < \dots < \alpha_j^{n_j} < 1$
\end{enumerate}
for all $x_j$ appearing in $D$. We write $E_\star$ to mean $(E, \{E_j^i\} ,\{\alpha_j^i\})$, and we call $D$ the parabolic divisor of $E_\star$. For a fixed $x_j$, we call $\{E_j^i\}$ the flag over $x_j$ and $\{\alpha_j^i\}$ the weights over $x_j$. We call $\alpha_j^i$ the weight associated to $E_j^i$.
\end{definition}

Given a parabolic bundle $E_\star$, we get induced parabolic structures on subquotients of $E_\star$.

\begin{definition}\label{def: definition_parabolic_subquotients}
Let $E_\star$ be a parabolic bundle on $X$ with parabolic divisor $D$, and let $F \subseteq E$ be a subbundle of $E$. The vector bundles $F$ and $E/F$ carry an induced parabolic structure with parabolic divisor $D$ as follows: 
\begin{enumerate}
	\item The flag of $F_{x_j}$ over $x_j$ is given by 
	\begin{align*}
	F_{x_j} = E_j^1 \cap F_{x_j} \supseteq E_j^2 \cap F_{x_j} \supseteq \dots \supseteq E_j^{n_j + 1} \cap F_{x_j} = 0
	\end{align*}
	after removing redundancies, and the weight associated to $F_j^i$ is given by
        \[\max_{1 \leq k \leq n_j} \left\{\alpha_j^k : F_j^i = E_j^k \cap F_{x_j} \right\}.\]
	\item The flag of $(E/F)_{x_j}$ is given by
	\begin{align*}
		(E/F)_{x_j} = \left(E^1_j + F_{x_j}\right)/F_{x_j} \supseteq \left(E^2_j + F_{x_j}\right)/F_{x_j} \supseteq \dots \supseteq \left(E^{n_j + 1}_j + F_{x_j}\right)/F_{x_j} = 0
	\end{align*}
	after removing redundancies, and the weight associated to $(E/F)_{j}^i$ is given by \[\max_{1 \leq k \leq n_j}\left\{\alpha_j^k : (E/F)_{j}^i = \left(E_j^k + F_{x_j}\right)/F_{x_j}\right\}.\]
\end{enumerate}
\end{definition}

\begin{definition}\label{def: definition_parabolic_degree}
Let $E_\star$ be a parabolic vector bundle on $X$ with parabolic divisor $D = x_1 + \dots + x_n$. The \emph{parabolic degree} is the quantity
\begin{align*}
	\pardeg (E_\star) &:= \deg(E) + \sum_{j = 1}^n \sum_{i = 1}^{n_j} \alpha_j^i \cdot \dim\left(E_j^i/E_j^{i + 1}\right), 
\end{align*}
and we call $\mu_\star(E_\star) := \pardeg (E_\star)/\rank (E)$  the \emph{parabolic slope} of $E_\star$.
\end{definition}

To define the notion of parabolic tensor product and morphisms between parabolic bundles, we will need the notion of a parabolic sheaf. To that end, let $\br$ be the category consisting of the real numbers with a single morphism $i^{a,b}\colon a \to b$ if $a \leq b$.

\begin{definition}\label{def: R-filtered_module}
	An \emph{$\br$-filtered $\mathscr{O}_X$-module} is a functor $E\colon \br \to \Mod_{\mathscr{O}_X}$, where $\Mod_{\mathscr{O}_X}$ is the category of $\mathscr{O}_X$-modules. We write $E_\alpha$ for $E(\alpha)$, and we write $E_\star$ for the functor $E$. We write $i_E^{a, b} = E(i^{a,b})$.
	
	Set $E[\alpha]_\star$ to be the functor given by $E[\alpha]_\beta = E_{\alpha + \beta}$, and $i^{a , b}_{E[\alpha]} = i^{a + \alpha, b + \alpha}_E$. We write $i^{[\alpha, \beta]}_E \colon E[\alpha]_\star \to E[\beta]_\star$ for the natural transformation so that $i_E^{[\alpha, \beta]}(\gamma) = i_E^{\alpha + \gamma, \beta + \gamma}$. If $f \colon E_\star \to F_\star$ is a natural transformation of $\br$-filtered $\mathscr{O}_X$-modules, we write $f[\alpha] \colon E[\alpha]_\star \to F[\alpha]_\star$ for the induced map.
\end{definition}

\begin{definition}\label{def: definition_parabolic_sheaf}
A \emph{parabolic sheaf} $E_\star$ with divisor $D$ is an $\br$-filtered $\mathscr{O}_X$-module with an isomorphism $j_E \colon E_\star (-D) \to E_\star[1]$ such that $i_E^{[0,1]} \circ j_E = \id_{E_\star} \otimes i_D$, where $i_D \colon \mathscr{O}_X(-D) \to \mathscr{O}_X$ is the inclusion.

A natural transformation $f\colon E_\star \to F_\star$ is a parabolic morphism if $f[1] \circ j_E = (f \otimes \id )\circ  j_F$. Let $\hom(E_\star, F_\star)$ be the set of parabolic morphisms and $\shom(E_\star, F_\star)$ be the sheaf of parabolic morphisms defined by $\shom(E_\star,F_\star)(U) = \hom({E_\star}|_U, {F_\star}|_U)$. The sheaf $\shom(E_\star, F_\star)$ has a parabolic structure given by $\shom(E_\star, F_\star)_\alpha = \shom(E_\star, F[\alpha]_\star)$ with morphism $\shom(E_\star, F_\star)_\alpha \to \shom(E_\star, F_\star)_\beta$ induced by $F[\alpha]_\star \to F[\beta]_\star$.
\end{definition}

The notions of parabolic vector bundle and parabolic sheaf are related.

\begin{example}\label{ex: parabolic_bundles_to_sheaves}
	Let $(E, \{E_j^i\}, \{\alpha_j^i\})$ be a parabolic vector bundle with parabolic divisor $D$. For $0 \leq \alpha < 1$, write $\beta(\alpha,j ) = \min\{i : \alpha_j^i \geq \alpha\}$ for all $\alpha \leq \max \{\alpha_j^i\}$ and $\beta(\alpha,j) = n_j + 1$ for $\max \{\alpha_j^i\} < \alpha < 1$. Then set 
	\begin{align*}
		E_\alpha &= \bigcap_{j = 1}^n \ker \left(E \to E_{x_j}/ E_j^{\beta(\alpha, j)}\right).
	\end{align*}
	We extend outside $[0,1)$ by $E_{x + n} = E_x(-nD)$ for all $n \in \bz$, with morphisms $i^{\alpha, \beta} \colon E_\alpha \to E_\beta$ given by their natural inclusions.
\end{example}

\begin{definition}\label{def: definition_parabolic_tensor}
Let $E_\star$ and $F_\star$ be two parabolic sheaves with the same parabolic divisor $D$, and let $j \colon C \setminus D \hookrightarrow C$ be the natural inclusion. Suppose that $E_\alpha$ and $F_\alpha$ are both locally free for all $\alpha$. Then, we define the \emph{parabolic tensor product} $E_\star \otimes F_\star$ to be the parabolic sheaf such that $(E_\star \otimes F_\star)_\alpha$ is generated by all $E_a \otimes F_b$ where $a + b = \alpha$ when viewed as a subbundle of $j_\star j^* (E \otimes F)$, and morphisms $i_{E_\star\otimes F_\star}^{\alpha, \beta}$ to be the natural inclusion maps.
\end{definition}

A parabolic sheaf $E_\star$ is constant along intervals $[\alpha_i, \alpha_{i + 1})$. To each parabolic sheaf, we can associate a coparabolic sheaf $\widehat{E}_\star$ that is constant along $(\alpha_i, \alpha_{i + 1}]$. (See \cite[Figure 1]{Boden_Yokogawa_moduli_of_parabolic_higgs}). 

\begin{definition}\label{def: definition_coparabolic_sheaf}
	Let $E_\star$ be a parabolic sheaf. We define $\widehat{E}_\star$ by the rule
	\begin{align*}
		\widehat{E}_x &= \begin{cases}
			E_x, & x \neq \alpha_i,\\
			E_{\alpha_{i + 1}}, & x = \alpha_i,
		\end{cases}
	\end{align*}
	and we call $\widehat{E}_\star$ the \emph{coparabolic sheaf} associated to $E_\star$. If $E_\star$ is a parabolic bundle, then we call $\widehat{E}_\star$ the \emph{coparabolic bundle} associated to $E_\star$.
\end{definition}

\begin{proposition}\label{prop: tensor_hom_isomorphism}
  If we write $F_\star^\vee$ to mean $\shom(F_\star, \mathscr{O}_X)_\star$ $($where $\mathscr{O}_X$ has the trivial parabolic structure$)$, then $\shom(E_\star, F_\star)_\star \cong (E_\star^\vee \otimes F_\star)_\star$.
\end{proposition}
\begin{proof}
	See \cite[Lemma 3.6]{Yokogawa_infinitesimal_deformation}.
\end{proof}

Let $E_\star$ and $F_\star$ be parabolic vector bundles with parabolic weights $\{\alpha_i^j\}$ and $\{\beta_i^j\}$. A morphism $f \colon E_\star \to F_\star$ is simply a morphism of vector bundles such that at every point $x_j \in D$, the morphism preserves parabolic structure. That is, if $\alpha_j^i \geq \beta_j^k$, then $f(E_j^i) \subseteq F_j^k$ .  A \emph{strongly} parabolic morphism is a parabolic morphism such that if $\alpha_j^i \geq \beta_j^k$, then $f(E_j^i) \subseteq F_j^{k + 1}$. In other words, strongly parabolic morphisms are sections of $\widehat{\shom(E_\star, F_\star)}$.

As an example, we spell out the case where $E_\star$ and $F_\star$ are parabolic line bundles. If $E_\star$ and $F_\star$ are rank~$1$ parabolic bundles, then $\shom(E_\star, F_\star) \cong E^\vee \otimes F (-D')$, where $D' \subseteq D$ is the divisor consisting of points $x_j$ where the weight of $E_\star$ is greater than the weight of $F_\star$. Similarly, $\widehat{\shom(E_\star, F_\star)} \cong E^\vee \otimes F(-D'')$, where $D''$ is the divisor of points where the weight of $E_\star$ is greater than or equal to the weight of $F_\star$.

\begin{definition}\label{def: definition_parabolic_higgs}
A \emph{parabolic Higgs bundle} $(E_\star, \theta)$ on $(X,D)$ is a parabolic bundle $E_\star$ on $(X,D)$ along with an $\End(E_\star)$-valued logarithmic $1$-form $\theta \in H^0(X, \End(E_\star) \otimes \Omega_X^1(\log D))$. The map $\theta$ is called the Higgs field. Here, ``logarithmic'' means $\theta$ has simple poles along the divisor $D$.

We say that $(E_\star,\theta)$ is a \emph{strongly parabolic Higgs bundle} if it is a parabolic Higgs bundle and $\theta$ is a strongly parabolic morphism $\theta\colon E_\star \to E_\star \otimes \Omega_X^1(\log D)$. In other words, $\theta$ is a global section of the bundle $\widehat{\End(E_\star)} \otimes \Omega_X^1(\log D)$.
\end{definition}

We consider strongly parabolic Higgs bundles to use deformation-theoretic techniques. These strongly parabolic Higgs bundles were also considered in the work of Tholozan--Toulisse \cite{Tholozan_Toulisse}.

There are notions of stability and semi-stability for parabolic Higgs bundles.

\begin{definition}\label{def: definition_higgs_parabolic_stable}
If $(E_\star, \theta)$ is a Higgs bundle with a parabolic structure, then we say that $(E_\star, \theta)$ is stable (resp.\ semi-stable) if all sub-Higgs bundles $(F_\star, \theta|_F)$ satisfy $\mu_\star(F_\star) < \mu_\star(E_\star)$ (resp.\ $\mu_\star(F_\star) \leq \mu_\star(E_\star)$). We say that $(E_\star, \theta)$ is \emph{poly-stable} if it is isomorphic to a direct sum of stable Higgs bundles.
\end{definition}

\begin{lemma}\label{lem: parabolic_tensor_degree}
If\, $F_\star$ and $G_\star$ are two parabolic bundles on $\bp^1$, then
\begin{align}
\pardeg \shom(E_\star, F_\star)_\star = \rank E \cdot \pardeg F_\star - \rank F \cdot \pardeg E_\star. 
\end{align}
\end{lemma}
\begin{proof}
	See \cite[Section 2]{Biswas_principal_bundles} for a discussion on computing $\pardeg (E_\star \otimes F_\star)_\star$ for parabolic bundles $E_\star$ and $F_\star$. Then to get the statement for $\shom(E_\star, F_\star)_\star$, we note that $\shom(E_\star, F_\star)_\star = ((E_\star)^\vee \otimes F_\star)_\star$.
\end{proof}

The parabolic degree satisfies properties similar to those satisfied by the usual notion of degree.

\begin{lemma}\label{lem: par_deg_properties}
  Let
  $$0 \lra F_\star^1 \overset{d_1}\lra F_\star^2 \overset{d_2}\lra \cdots \xrightarrow{d_{N-2}} F_\star^{N - 1} \xrightarrow{d_{N-1}}F_\star^N \lra 0$$ be an exact sequence of parabolic bundles. Then,
\begin{align*}
	\sum_{k = 1}^N (-1)^k \pardeg F_\star^k = 0.
\end{align*}
In particular, if
$$0 \lra F_\star \lra E_\star \lra (E/F)_\star \lra 0$$ is a short exact sequence of parabolic bundles, then $\pardeg(E_\star) = \pardeg(F_\star) + \pardeg((E/F)_\star)$. A parabolic bundle $E_\star$ is semi-stable $($resp.\ stable$)$ if and only if for every quotient $E_\star \to (E/F)_\star$,  $\mu_\star(E_\star) \leq \mu_\star((E/F)_\star)$ $($resp.\ $\mu_\star(E_\star)< \mu_\star((E/F)_\star))$.
\end{lemma}
\begin{proof}
	See \cite[Lemma 2.4.5]{Landesman-Litt-very_general_curves} for the claim about short exact sequences. We deduce the statement about general exact sequences from the claim about short exact sequences. 
	
	Let 
	\begin{align*}0 \overset{d_0}\lra F_\star^1 \overset{d_1}\lra F_\star^2 \overset{d_2}\lra \cdots \xrightarrow{d_{N-2}} F_\star^{N - 1} \xrightarrow{d_{N-1}}F_\star^N \overset{d_{N}}\lra 0	
	\end{align*}
 	be an exact sequence of parabolic bundles. Then as parabolic bundles, $\ker(d_{k+1})_\star = \im(d_{k})_\star$. We now break up the exact sequence into many short exact sequences. We have exact sequences of parabolic bundles
 	\begin{align*}
 		0 \lra \ker(d_{k})_\star \lra F_\star^{k} \lra \im(d_{k})_\star \lra 0.
 	\end{align*}
 	Then since we know that parabolic degrees are additive in short exact sequences, we have that
 	\begin{align*}
 		\pardeg F^{k}_\star = \pardeg \ker(d_{k})_\star + \pardeg \im(d_{k})_\star.
 	\end{align*}
 	Then,
 	\begin{align*}
 		\sum_{k = 1}^N (-1)^k \pardeg F_\star^k &= \sum_{k = 1}^N (-1)^k \left (\pardeg \ker(d_{k})_\star + \pardeg \im(d_{k})_\star \right )\\
 		&= \sum_{k = 1}^N (-1)^k \left (\pardeg \im(d_{k-1})_\star + \pardeg \im(d_{k})_\star \right ).
 	\end{align*}
 	This is a telescoping sum, and so simplifying yields
 	\begin{align*}
\sum_{k = 1}^N (-1)^k \left (\pardeg \im(d_{k-1})_\star + \pardeg \im(d_{k})_\star \right ) = - \pardeg \im(d_0)_\star + (-1)^N\pardeg \im(d_{N})_\star.
 	\end{align*}
 	But $d_0 = d_N = 0$, so $- \pardeg \im(d_0)_\star + (-1)^N\pardeg \im(d_{N})_\star =0$, proving the claim.
\end{proof}

\begin{lemma}\label{lem: bounds_pardeg}
Let $E_\star$ be a parabolic vector bundle with parabolic divisor $D = x_1 + \dots + x_n$. Then
\begin{align*}
	\deg E\leq \pardeg E_\star < \deg E + \rank E \cdot n.
\end{align*}
\end{lemma}
\begin{proof}
	By \Cref{def: definition_parabolic_degree}, 
	\begin{align*}
		\pardeg E_\star = \deg E + \sum_{j = 1}^n \sum_{i = 1}^{n_j} \alpha_i^j \cdot \dim\left(E_j^i/E_{j}^{i + 1}\right).
	\end{align*}
	Since our parabolic weights $\alpha_i^j$ lie in $[0,1)$, 
	\begin{align*}
		0 \leq \sum_{j = 1}^n \sum_{i = 1}^{n_j} \alpha_i^j \cdot \dim\left(E_j^i/E_{j}^{i + 1}\right) & < \sum_{j = 1}^n \sum_{i = 1}^{n_j} \dim\left(E_j^i/E_{j}^{i + 1}\right) = \sum_{j = 1}^n \dim\left(E_{x_j}\right). 
	\end{align*}
	But $E$ is a vector bundle, so the fibers have dimension equal to $\rank E$. Therefore, 
\begin{equation*}\pushQED{\qed}	
		\deg E \leq \pardeg E_\star < \deg E + \sum_{j = 1}^n \dim \left(E_{x_j}\right) = \deg V +\rank E \cdot n.\qedhere \popQED
	\end{equation*}
\renewcommand{\qed}{}     
\end{proof}

We seek to study parabolic Higgs bundles on $(X,D)$ because they correspond to local systems on $X\setminus D$. This correspondence is due to the work of several people including Corlette, Donaldson, Hitchin, and Simpson. For our purposes, we use nonabelian Hodge theory as described in \cite{Simpson_harmonic_bundles}, and we include it as \Cref{prop: Non-Abelian_Hodge_Theorem} for convenience. We will refer to \Cref{prop: Non-Abelian_Hodge_Theorem} for the correspondence.

\begin{theorem}\label{prop: Non-Abelian_Hodge_Theorem}
Let $X$ be a smooth proper curve and $D$  a reduced effective divisor. There is a one-to-one correspondence between poly-stable parabolic Higgs bundles of parabolic degree zero on $X$ with parabolic divisor $D$ and the semisimple local systems on $X\setminus D$. Under this correspondence, the poly-stable parabolic Higgs bundles with zero Higgs field correspond to the unitary local systems on $X$, and strongly parabolic Higgs bundles correspond to the local systems with unitary local monodromy at the punctures of $X \setminus D$.
\end{theorem}
\begin{proof}
	See \cite[p.~718, Theorem]{Simpson_harmonic_bundles}.
\end{proof}

\begin{remark}
	The local systems with unitary monodromy must have unitary local monodromy at the punctures, but the converse is not true. In the context of \Cref{prop: Non-Abelian_Hodge_Theorem}, this is captured by the phenomenon that there are strongly parabolic Higgs bundles with nonzero Higgs field. These strongly parabolic Higgs bundles with nonzero Higgs field correspond to local systems with unitary monodromy locally at the punctures but are not globally unitary.
\end{remark}

We will describe which Higgs bundles correspond to local systems underlying a complex variation of Hodge structure in \Cref{sec:VHS}.

\section{Variations of Hodge structure}\label{sec:VHS}

We give background on variations of Hodge structure and their associated Higgs bundles in this section. We primarily use the notation of \cite[Section 4]{Landesman-Litt-very_general_curves}.

\begin{definition}\label{def: definition_VHS}
Let $X$ be a smooth irreducible variety over $\bc$ and $w \in \bz$ an integer. A \emph{complex variation of Hodge structure} ($\bc$-VHS) of \emph{weight} $w$ is a triple $(V, V^{p,q}, \overline{\partial})$, where $V$ is a complex $C^\infty$-bundle on $X$, $\oplus_{p + q = w}V^{p,q} = V$ is a direct sum decomposition, and $\overline{\partial}$ is a flat connection satisfying Griffiths transversality. That is,
\begin{align*}
	\overline{\partial}\left(V^{p,q}\right) \subseteq A^{1,0}(V^{p,q}) \oplus A^{0,1}\left(V^{p,q}\right) \oplus A^{1,0}\left(V^{p - 1, q + 1}\right) \oplus A^{0,1}\left(V^{p + 1, q - 1}\right).
\end{align*}
Here, $A^{i,j}(V^{p,q})$ is the sheaf of $V^{p,q}$-valued $(i,j)$-forms. A \emph{polarization} on $(V, V^{p,q}, \overline{\partial})$ is a flat Hermitian form $\psi$ on $V$ such that the $V^{p,q}$ are orthogonal to each other under $\psi$ and that on each $V^{p,q}$, the form $(-1)^p \psi$ is positive definite. A \emph{polarizable} $\bc$-VHS is a $\bc$-VHS that admits a polarization $\psi$.

Given a $\bc$-VHS $(V, V^{p,q}, \overline{\partial})$, we call the holomorphic flat bundle $(E, \nabla):= (\ker (\overline{\partial}) \otimes_\bc \mathscr{O}_X, \id \otimes d)$ the \emph{holomorphic flat bundle} associated to the $\bc$-VHS. We get a filtration $F^p E$ on $E$ induced by the filtration $F^p V = \oplus_{j \geq p} V^{j, p + q - j}$ on $V$, known as the Hodge filtration. The filtration satisfies Griffiths transversality:
\begin{align*}
	\nabla(F^p) \subseteq F^{p - 1} \otimes \Omega_X^1(\log D).
\end{align*}
If $\mathbb{V}$ is a local system isomorphic to $\ker (\overline{\partial})$ for some polarizable $\bc$-VHS $(V, V^{p,q}, \overline{\partial})$, then we say that $\mathbb{V}$ underlies a polarizable $\bc$-VHS.
\end{definition}

\begin{remark}
	To an irreducible local system $\bv$ underlying a $\bc$-VHS,  there is a $\bc$-VHS underlying $\bv$ of weight~$w$ for every $w \in \bz$. To see this, we pick some $w \in \bz$, and we let $(V, V^{p,w - p}, \overline{\partial})$ be a $\bc$-VHS underlying $\bv$ of weight $w$. For any $k \in \bz$, we let  $(V', (V')^{p,w + k - p} , \overline{\partial '})$ be a $\bc$-VHS, where $V' = V$, $(V')^{p,w + k-p} = V^{p,w - p}$, and $\overline{\partial '} = \overline{\partial}$. Then, the weight of $(V', (V')^{p,w + k - p} , \overline{\partial '})$ is $w + k$, and Griffiths transversality still holds. Since the local system determined by $(V', (V')^{p,-p} , \overline{\partial '})$ is given by $\ker (\overline{\partial '})$ and $\overline{\partial '} = \overline{\partial}$, the two $\bc$-VHS underlie the same local system.
	
	The data of a weight $w$ is extra data. Often this choice of weight is irrelevant data, and we will often ask for a local system to underlie a $\bc$-VHS without specifying the weight because the weight will not matter.
\end{remark}

\begin{definition}\label{def: definition_Deligne_canonical_extension}
Let $X = \overline{X}\setminus D$ with $D$ a simple normal crossings divisor, and let $(E,\nabla)$ be a flat holomorphic vector bundle on $X$. The \emph{Deligne canonical extension} of $(E,\nabla)$ is the unique logarithmic flat bundle $(\overline{E}, \overline{\nabla})$ on $\overline{X}$ with regular singularities along $D$, along with an isomorphism $(\overline{E}, \overline{\nabla})|_X \cong (E, \nabla)$, and such that the eigenvalues of the residues of $\overline{\nabla}$ along $D$ have real part in $[0,1)$.
\end{definition}

The existence of the Deligne canonical extension is a theorem of Deligne; see \cite[Remarques 5.5(i)]{Deligne_canonical_extension}.

\begin{definition}\label{def: definition_parabolic_bundle_to_Deligne_canonical_extension}
Let $X = \overline{X} \setminus D$ be a curve, where $\overline{X}$ is some smooth proper curve and $D$ is a reduced effective divisor. Let $(E, \nabla)$ be a flat holomorphic vector bundle of rank $n$ on $X$ with regular singularities along $D$. Suppose at some $x_j$, the eigenvalues of the residue matrix $ \res(\nabla)(x_j)$ are given by $\eta_j^1,\dots, \eta_j^n$. Then, the real parts $\Re(\eta_j^i)$ lie in $[0,1)$. We write $ \alpha_j^i = \Re(\eta_j^i)$ and order the $\alpha_j^i$ in increasing order and remove repetitions so
\begin{align*}
	0 \leq \alpha_j^1 < \dots < \alpha_j^{n_j} < 1.
\end{align*}
Let $E_j^i \subseteq E_{x_j}$ be the sum of all of the generalized eigenspaces of $\res(\nabla)(x_j)$ such that the real parts of the associated eigenvalues are at least $\alpha_j^i$. Then, $E = (E, \{ E_j^i\}, \{\alpha_j^i\})$ is the \emph{parabolic bundle associated to the connection} $\nabla$.
\end{definition}

\begin{definition}\label{def: definition_associated_higgs_bundle}
Let $\overline{X}$ be a smooth projective curve, $D$ a reduced effective divisor, and $X = \overline{X} \setminus D$. Let $(V, V^{p,q}, \overline{\partial})$ be a $\bc$-VHS on $X$ and $(E, \nabla)$ its associated holomorphic flat bundle with Hodge filtration $F^\bullet$. Let $(\overline{E}_\star, \overline{\nabla})$ be its associated parabolic bundle on $X$ in the sense of \Cref{def: definition_Deligne_canonical_extension,def: definition_parabolic_bundle_to_Deligne_canonical_extension}. Then by  \cite[Section 7]{Brunebarbe}, there is a canonical extension of the Hodge filtration $F^\bullet$ to $(\overline{E}, \overline{\nabla})$. Let $\overline{F}^\bullet$ be this extension.

We call $(\oplus_p \gr_{\overline{F}^\bullet}^p \overline{E}_\star, \oplus \gr_{\overline{F}^\bullet}^p \overline{\nabla})$ the \emph{parabolic Higgs bundle associated} to the $\bc$-VHS. The $\mathscr{O}_X$-linearity of the map $\oplus \gr_{\overline{F}^\bullet}^p \overline{\nabla}$ is due to Griffiths transversality of $\overline{\nabla}$ and $\overline{F}^\bullet$.
\end{definition}

\begin{theorem}\label{prop: non-abelian Hodge theorem for C-VHS}
Let $X = \overline{X} \setminus D$, where $\overline{X}$ is a smooth proper curve and $D$ is a reduced effective divisor. Let $(V, V^{p,q}, \overline{\partial})$ be a complex variation of Hodge structure on $X$, $(E, F^\bullet, \nabla)$ be its associated holomorphic flat bundle with Hodge filtration $F^\bullet$, and $(\oplus E^p_\star, \oplus \theta_p) :=(\oplus_p \gr_{\overline{F}^\bullet}^p \overline{E}_\star, \oplus \gr_{\overline{F}^\bullet}^p \overline{\nabla})$ be its associated parabolic Higgs bundle defined in \Cref{def: definition_associated_higgs_bundle}. Then,
\begin{enumerate}
	\item $\oplus E_\star^p$ is of parabolic degree zero;
	\item $(\oplus E_\star^p, \oplus \theta_p)$ is parabolic semi-stable as a Higgs bundle;
	\item if the local system $\ker(\nabla)$ is irreducible, then $(\oplus E_\star^p, \oplus \theta_p)$ is parabolic stable as a Higgs bundle.
\end{enumerate}
\end{theorem}
\begin{proof}
	This is due to Simpson \cite[Theorem 5]{Simpson_harmonic_bundles}.
\end{proof}

\Cref{prop: non-abelian Hodge theorem for C-VHS} produces a parabolic Higgs bundle from a $\bc$-VHS. The parabolic Higgs bundle is the same as the one given by the correspondence given in \Cref{prop: Non-Abelian_Hodge_Theorem}.

\begin{definition}\label{def: definition_generic_distinct}
Let $X$ be a smooth projective curve and $D = x_1 + \dots + x_n$ a reduced effective divisor on~$X$. We fix a positive integer $r$. Let $\{\alpha_j^i\}$ be a collection of parabolic weights. 
\begin{enumerate}
	\item\label{d:dgd-1} We say that $\{\alpha_j^i\}$ is \emph{smooth} with respect to $(X,D,r)$ if every semi-stable parabolic Higgs bundle of rank $r$ on $(X,D)$ of parabolic degree zero with  weights $\{\alpha_j^i\}$ is actually stable.
	\item\label{d:dgd-2} We say that the collection of weights $\{\alpha_j^i\}$ is \emph{distinct} with respect to $(X,D,r)$ if for all $j$, we have that $n_j = r$, where $n_j$ is the number of weights over $x_j$.
	\item\label{d:dgd-3} We say that the weights $\{\alpha_j^i\}$ are \emph{generic} with respect to $(X,D,r)$ if \eqref{d:dgd-1} and \eqref{d:dgd-2} are both satisfied.
\end{enumerate}
Let $\mathbb{V}$ be a local system on $X \setminus D$ with unitary local monodromy around the punctures. We say that the local monodromy of $\mathbb{V}$ is generic if the collection of associated weights of the parabolic Higgs bundle is generic (via nonabelian Hodge theory).

Given some local monodromy data $C = (C_1,\dots, C_n)$, we say that the $n$-tuple of conjugacy classes $C$ is \emph{connected} if the following holds. For any subspace $V \subseteq \bc^{r}$ and any collection of matrices $(A_1,\dots, A_n)$ with $A_i \in C_i$ such that $V$ is stable under all of the matrices $A_i$ and
\begin{align*}
	\prod_{i = 1}^n \det (A_i|_V) = 1,
\end{align*}
we have that either $V = 0$ or $V = \bc^r$.
\end{definition}

We call an $n$-tuple of semisimple conjugacy classes $C = (C_1,\dots, C_n)$ connected because this condition ensures that the relative character variety $\chi(\pi_1(\Sigma_{g,n}), \gl_r(\bc), C)$ (when the $C_i$ are semisimple) is connected by \cite[Theorem 5.1.1]{Hausel_Letellier_Rodriguez-Villegas}. We hope that this nonstandard naming convention is not too confusing for the reader.

\begin{example}\label{ex: generic weights example}
  Let $(X,D,r)$ be fixed. Let $\{\alpha_j^i\} \subseteq [0,1)$ be a set of size $\deg D \cdot r$ such that no proper subset of  $\{\alpha_j^i\}$ has elements that sum to an integer. Then, $\{\alpha_j^i\}$ is generic with respect to $(X,D,r)$ in the sense of \Cref{def: definition_generic_distinct}\eqref{d:dgd-3}.
    To check this, we take any parabolic semi-stable Higgs bundle $(E_\star, \theta)$ of parabolic degree zero and observe that any sub-Higgs bundle $(F_\star, \theta|_F)$ of $(E_\star, \theta)$ satisfies
\begin{align*}
	\frac{\pardeg F_\star}{\rank F} = \mu_\star(F_\star) \leq \mu_\star(E_\star) = \frac{\pardeg E_\star}{\rank E} = 0.
\end{align*}
By definition, the parabolic degree of $F_\star$ is given by
\begin{align}
	\pardeg F_\star = \deg F + \sum_{j = 1}^{\deg D} \sum_{k = 1}^{\rank F}\alpha_j^{{i_k}}.
\end{align}
Here,  $\{ \alpha_j^{i_k}\}$ is the collection of weights of $F_\star$. We note that
\begin{align*}
	\deg F + \sum_{j = 1}^{\deg D} \sum_{k = 1}^{\rank F}\alpha_j^{{i_k}} < 0
\end{align*}
because by assumption, the expression
\begin{align*}
	\sum_{j = 1}^{\deg D} \sum_{k = 1}^{\rank F}\alpha_j^{{i_k}}
\end{align*}
is not an integer and $\deg F$ is always an integer.
\end{example}

\begin{remark}
	In terms of the local monodromy, \Cref{def: definition_generic_distinct}(\ref{d:dgd-2}) says that at every point, the eigenvalues are all distinct. Part~(\ref{d:dgd-1}) guarantees that the local systems in $\chi(\pi_1(X \setminus D), \gl_r(\bc), C)$ are irreducible.
\end{remark}

Let $X$ be a curve. Given a parabolic Higgs bundle $(E_\star, \theta)$ on $(X,D)$, the $\supth{i}$ Higgs cohomology group of $(E_\star, \theta)$ is defined to be the $\supth{i}$ hypercohomology of the complex
$$E_\star[0] \overset{\theta}\lra \widehat{E_\star} \otimes \Omega_X^1(\log D)[-1].$$
That is, we define the $\supth{i}$ Higgs cohomology to be
\begin{align*}
	\mathbb{H}^i\left (E_\star\xrightarrow{\theta } \widehat{E_\star} \otimes \Omega_X^1(\log D)\right).
\end{align*}
We will often omit the cohomological degree that each term of the complex sits in if it is clear from the context. 

If $(E_\star, \theta)$ underlies a $\bc$-VHS of weight $w$, $E_\star = \oplus_{p + q = w} E^{p,q}_\star$, then the Higgs cohomology has a Hodge structure of weight $w + i$ with $(p,q)$ part given by
\begin{align*}
	\mathbb{H}^{p,q}\left (E_\star \xrightarrow{\theta } \widehat{E_\star} \otimes \Omega_X^1(\log D) \right) &:= \mathbb{H}^i\left (E^{p,q - i}_\star\xrightarrow{\theta }\widehat{E^{p - 1, q - i + 1}_\star} \otimes \Omega_X^1(\log D) \right).
\end{align*}
This means that
\begin{align*}
	\mathbb{H}^i \left (E_\star \xrightarrow{\theta } \widehat{E_\star} \otimes \Omega_X^1(\log D) \right) = \bigoplus_{p +q = w + i} \mathbb{H}^i\left (E^{p,q - i}\xrightarrow{\theta }\widehat{E^{p - 1, q - i + 1}_\star} \otimes \Omega_X^1(\log D) \right).
\end{align*}
For a review of hypercohomology, we suggest \cite[Example A.28]{Peters_Steenbrink}.

\begin{definition}\label{def: definition_Hodge_length}
Let $(E_\star, \theta)$ be a parabolic Higgs bundle underlying a $\bc$-VHS on a curve $X$ with parabolic divisor $D$. We say that $\mathbb{H}^1(E_\star \xrightarrow{\theta} \widehat{E_\star} \otimes \Omega_X^1(\log D))$ is of \emph{Hodge length} $\ell$ if $\mathbb{H}^{k, 1 - k}(E_\star \xrightarrow{\theta} \widehat{E_\star} \otimes \Omega_X^1(\log D)) = 0$ for all $k < - \ell + 1$ and $k > \ell$.
\end{definition}

\section{Minimal energy local systems}\label{sec:minimal energy}

Let $X$ be a smooth proper curve and $D = x_1 + \dots + x_n$ a reduced effective divisor on $X$ (\textit{i.e.}, a collection of points). We let $M = M(X, D,r, \{\alpha_i^j\})$ be the coarse moduli space of semi-stable strongly parabolic Higgs bundles of rank $r$ of parabolic degree zero on $X$ with parabolic divisor $D$ and parabolic weights $\{\alpha_i^j\}$ (\Cref{def: definition_parabolic_higgs,def: definition_higgs_parabolic_stable}). This space is constructed in \cite[Section 2]{Yokogawa_compact_moduli_space}. As $(X,D,r)$ are fixed, any smoothness, distinctness, and genericity conditions on a collection of weights $\{\alpha_j^i\}$ are assumed to be with respect to $(X,D,r)$ in the sense of \Cref{def: definition_generic_distinct}.

As before, for a semi-stable parabolic Higgs bundle $(E_\star, \theta)$ to be in $M$, we require that the Higgs field $\theta$ is a strongly parabolic Higgs field. This ensures that when we deform a Higgs bundle (in $M$) corresponding to a local system with unitary local monodromy, the resulting Higgs bundle gives us a local system with unitary local monodromy. 

To see why this is necessary, we consider an example. For every $t \in \bc$, consider a local system $\bv_t$ on $\Sigma_{g,n}$ with local monodromy at the puncture $x_1$ given by the matrix
\begin{align*}
	\begin{pmatrix}
		1 & t \\ 0 & 1
	\end{pmatrix}.
\end{align*}
Whenever $t \neq 0$, the local monodromy stays in the same conjugacy class in $\sl_2(\br)$. However,  the limit $\bv_t$ does not exist in the relative character variety because $\bv_0$ has trivial monodromy at $x_1$. Taking unitary local monodromy (and therefore working with strongly parabolic Higgs bundles) fixes this because the semisimple conjugacy classes are closed.

\begin{remark}\label{rem: definition_adjoint}
Let $(E_\star, \theta) = (\oplus E^p, \oplus \theta_p)$ be stable parabolic Higgs bundle underlying an (irreducible) $\bc$-VHS. The adjoint Higgs bundle $(\ad(E_\star), \ad(\theta))$ is given by the quotient $\End(E_\star)/ \mathscr{O}_X$. Here, we are viewing $\mathscr{O}_X$ as the subsheaf of $\End(E_\star)$ of scalar endomorphisms of $E_\star$ with its trivial parabolic structure.

The adjoint Higgs bundle $\ad(E_\star)$ also underlies a $\bc$-VHS (of weight zero). This has a canonical grading, where the $\supth{k}$ graded piece $\ad(E_\star)^{k,-k}$ is given by $\oplus \shom(E^p_\star, E^{p + k}_\star)_\star$. Then,
\begin{align*}
	\mathbb{H}^1 = \mathbb{H}^1( \ad(E_\star)\xrightarrow{\ad(\theta) } \widehat{\ad(E_\star)} \otimes \Omega_X^1(\log D))
\end{align*}
has a Hodge structure of weight $1$. In this case $\dim_\bc \mathbb{H}^{p,1-p} = \dim_\bc \mathbb{H}^{1-p, p}$ because $\ad(E_\star)$ is self-dual, so $\mathbb{H}^1$ has a real Hodge structure. This is because $\ad(E_\star)$ is polarized, so it is dual to its conjugate. But it is self-dual, so its conjugate is itself.

We note that instead of working with the adjoint Higgs bundle, one can instead work with $\End(E_\star)$. But for our purposes, this does not matter because we work with $\mathbb{H}^1$, and their cohomologies agree:
\begin{align*}
		\mathbb{H}^1\left( \ad(E_\star)\xrightarrow{\ad(\theta) } \widehat{\ad(E_\star)} \otimes \Omega_X^1(\log D)\right) = \mathbb{H}^1\left( \End(E_\star)\xrightarrow{ \ad(\theta)} \widehat{\End(E_\star)} \otimes \Omega_X^1(\log D)\right).
\end{align*}
\end{remark}

\begin{proposition}\label{rem: tangent_space_adjoint}
	If\, $(E_\star, \theta)$ is a smooth point of $M$, then the tangent space at $(E_\star, \theta)$ is isomorphic to $\mathbb{H}^1( \ad(E_\star) \to \widehat{\ad(E_\star)} \otimes \Omega_X^1(\log D))$. If\, $(E_\star, \theta)$ is stable, then it is a smooth point in $M$.
\end{proposition}
	\begin{proof}
		See \cite[Section 3]{Thaddeus_variation_moduli_parabolic_higgs} for the smoothness claim. A detailed proof that the tangent space at $(E_\star, \theta)$ is isomorphic to $\mathbb{H}^1( \ad(E_\star) \to \widehat{\ad(E_\star)} \otimes \Omega_X^1(\log D))$ to can be found in \cite[Proposition 5.2.1]{Bottacin_symplectic}.
	\end{proof}

We now define the notion of minimal energy as a condition on $\mathbb{H}^1( \ad(E_\star) \to \widehat{\ad(E_\star)} \otimes \Omega_X^1(\log D))$. By \Cref{rem: tangent_space_adjoint}, such a notion is actually a condition on the tangent space of $M$ at $(E_\star, \theta)$.

\begin{definition}\label{def: definition_minimal_energy}
Let $(E_\star, \theta)$ be a stable parabolic Higgs bundle on a smooth proper curve $X$ with parabolic divisor $D$ and parabolic weights $\{\alpha_i^j\}$ corresponding to a complex variation of Hodge structure. We say that $(E_\star, \theta)$ is a \emph{minimal energy} parabolic Higgs bundle (or is of minimal energy) if the Higgs cohomology group
\begin{align*}
	\mathbb{H}^1\left( \ad(E_\star) \to \widehat{\ad(E_\star)}\otimes \Omega_X^1(\log D)\right)
\end{align*}
of $M$ at $(E_\star, \theta)$ is of Hodge length $1$.

If $(V, V^{p,q}, \overline{\partial})$ is the $\bc$-VHS corresponding to the minimal energy Higgs bundle $(E_\star, \theta)$, then we say that $(V, V^{p,q}, \overline{\partial})$ is a \emph{minimal energy} variation of Hodge structure.

Similarly, if $\mathbb{V}$ is a local system on $X \setminus D$ with unitary local monodromy such that the corresponding parabolic Higgs bundle $(E_\star, \theta)$ (under nonabelian Hodge theory) is of minimal energy, then we say $\mathbb{V}$ is a \emph{minimal energy} local system (or is of minimal energy).
\end{definition}

\begin{remark}
We define minimal energy Higgs bundles $(E_\star, \theta)$ only when $(E_\star, \theta)$ is stable because $\ad(E_\star)$ has a canonical grading (and is a $\bc$-VHS of weight zero) when $(E_\star, \theta)$ is stable. If $(E_\star, \theta)$ is not stable, then in general $\ad(E_\star)$ will still be a $\bc$-VHS, but it will not have a canonical grading.
\end{remark}

By \Cref{rem: tangent_space_adjoint}, if our parabolic weights are chosen to be smooth (in the sense of \Cref{def: definition_generic_distinct}), then every semi-stable parabolic Higgs bundle with those weights is actually stable. In this case, our moduli space~$M$ is smooth.

\begin{remark}\label{rem: etymology_minimal_energy}
	In the nonparabolic setting, there is an energy functional on the moduli space of Higgs bundles given by $||(E, \theta)|| = i \int_X \tr(\theta \wedge \theta^*)$, where $\theta^*$ is the adjoint to $\theta$ with respect to a harmonic metric. If $(E, \theta)$ is stable, then this metric is unique. Collier and Wentworth \cite[Section 4.1]{Collier_Wentworth_conformal_limits} showed that stable Higgs bundles that minimize $||(E, \theta)||$ are those whose tangent space is of Hodge length $1$ by exploiting the fact that the energy functional is a Morse--Bott function on the moduli space of Higgs bundles; see \cite[Section 8]{Hitchin_self_duality}. The tangent space to the moduli space of Higgs bundles splits into a direct sum $\oplus_{j \in \bz} T_j$, where $T_j = \mathbb{H}^{j, 1-j}$; see \cite[Section 8]{Hitchin_Lie_groups_teichmuller}. On each piece $T_j$, the Hessian of $||\cdot||$ has eigenvalue $j$, and hence a Higgs bundle cannot minimize the energy functional unless $T_j = 0$ for $j < 0$. Symmetry of the decomposition of the tangent space implies that the only nonzero pieces are $T_0$ and $T_1$, which is precisely the Hodge length~$1$ condition described in \Cref{def: definition_minimal_energy}.
\end{remark}

There are some classes of local systems that are always of minimal energy.

\begin{example}\label{ex: example_unitary_minimal_energy}
A unitary local system $\mathbb{V}$ corresponds to a parabolic Higgs bundle underlying a $\bc$-VHS coming in one piece $E_\star = E^1_\star$ and Higgs field $\theta = 0$. This is a consequence of the Mehta--Seshadri correspondence for irreducible unitary local systems and stable parabolic bundles; see \cite[Theorem 4.1]{Mehta_Seshadri_moduli_of_bundles}. Then, $\mathbb{V}$ is a minimal energy local system.
\end{example}

\begin{example}\label{ex: example_rigid_minimal_energy}
An irreducible rigid local system is always of minimal energy. This is because rigid local systems are isolated points, and so they have zero tangent space.

We give an example of a rigid local system. A result of Katz \cite[Theorem 1.1.2]{Katz_rigid_local_systems} says that for a rank $r$ irreducible local system $\mathbb{V}$ on $\bp^1 \setminus \{p_1,\dots, p_n\}$ with local monodromy data (conjugate to) $(A_1,\dots, A_n)$, $\mathbb{V}$ is rigid if and only if
\begin{align*}
	(2 - n)r^2 + \sum_{i = 1}^n \dim Z(A_i) = 2, 
\end{align*}
where $Z(A_i)$ is the centralizer of $A_i$.

In the case where $n = 3$ and $r = 2$, Katz's condition requires that 
\begin{align*}
	\dim Z(A_1) + \dim Z(A_2) + \dim Z(A_3) = 6.
\end{align*}
For a $2 \times 2$ nonscalar matrix $B$, $\dim Z(B) = 2$. Hence, any irreducible rank $2$ local system on $\bp^1 \setminus \{0, 1, \infty\}$ with nonscalar local monodromy around the punctures must be rigid and therefore is of minimal energy.
\end{example}

There is a $\bc^\times$-action on $M$. If $t \in \bc^\times$ and $(E_\star, \theta) \in M$, then $t$ acts on $(E_\star, \theta)$ by $t \cdot (E_\star, \theta) = (E_\star, t \theta)$. Simpson \cite[Theorem 8]{Simpson_harmonic_bundles} shows that the fixed points of this action on $M$ are precisely the locus of parabolic Higgs bundles coming from a $\bc$-VHS under the correspondence in \Cref{prop: Non-Abelian_Hodge_Theorem}.

The $\bc^\times$-action on $M$ induces a $\bc^\times$-action the tangent space $\mathbb{H}^1( \ad(E_\star) \to \widehat{\ad(E_\star)}\otimes \Omega_X^1(\log D))$ at a $\bc^\times$-fixed point in $M$, and this $\bc^\times$-action interacts with the Hodge structure in a natural way.

\begin{proposition}\label{prop: action_on_tangent_space}
	The $\bc^\times$-action on $M$ induces a $\bc^\times$-action on the tangent space of $M$ at a fixed point $(E_\star, \theta)$. This induces a $\bc^\times$-action $($through the isomorphism in \Cref{rem: tangent_space_adjoint}\,$)$ on $\mathbb{H}^1( \ad(E_\star) \to \widehat{\ad(E_\star)} \otimes \Omega_X^1(\log D))$. The action of\, $t \in \bc^\times$ on each piece of the Hodge structure $\mathbb{H}^{p,1-p}(\ad(E_\star) \to \widehat{ \ad(E_\star)} \otimes \Omega_X^1(\log D))$ is given by $t^p$.
\end{proposition}
	\begin{proof}
		Let $(E_\star, \theta)$ be a fixed point. By \cite[Theorem 8]{Simpson_harmonic_bundles}, $(E_\star, \theta)$ comes from a $\bc$-VHS and hence has a graded Higgs bundle $(E_\star, \theta) = (\oplus E^p_\star, \oplus \theta^p)$. The grading is unique up to shifting. Once we pick a grading, we get a map $f\colon (\oplus E^p_\star, \oplus t\theta^p) \to (\oplus E^p_\star, \oplus \theta^p)$ given by $f|_{E^p} = t^{p}$. This is an isomorphism of parabolic Higgs bundles.
		
		The map $f$ naturally gives us a morphism of complexes
\begin{center}
\begin{tikzcd}
	{\ad(E_\star)} & {\widehat{\ad(E_\star)} \otimes \Omega_X^1(\log D)} \\
	{\ad(E_\star)} & {\widehat{\ad(E_\star)} \otimes \Omega_X^1(\log D)\rlap{,}}
	\arrow["{{t\ad(\theta)}}", from=1-1, to=1-2]
	\arrow["{{\tilde{f}}}"', from=1-1, to=2-1]
	\arrow["{{\tilde{f} \otimes \id}}", from=1-2, to=2-2]
	\arrow["\ad(\theta)"', from=2-1, to=2-2]
\end{tikzcd}
\end{center}
	where $\ad(\theta)$ is the induced Higgs field on $\ad(E_\star)$. Here, $\tilde{f}|_{\ad(E_\star)^{p, -p}} $ is given by multiplication by $t^{p}$. Hence, we get a map on hypercohomology
	\begin{align*}
		\mathbb{H}^1(\ad(E_\star) \lra \widehat{\ad(E_\star)} \otimes \Omega_X^1(\log D)) \lra \mathbb{H}^1( \ad(E_\star) \lra \widehat{\ad(E_\star)} \otimes \Omega_X^1(\log D))
	\end{align*}
	induced by $\tilde{f}$.
	
	In order to compute the $\bc^\times$-action on $\mathbb{H}^{p, 1-p}$, we consider the morphism of complexes
	\begin{center}
		\begin{tikzcd}
	{\ad(E_\star)^{p,-p}} & {\widehat{\ad(E_\star)^{p-1,1-p}} \otimes \Omega_X^1(\log D)} \\
	{\ad(E_\star)^{p,-p}} & {\widehat{\ad(E_\star)^{p-1, 1-p}} \otimes \Omega_X^1(\log D)}\rlap{.}
	\arrow["{t\ad(\theta)}", from=1-1, to=1-2]
	\arrow["{\cdot t^{p}}"', from=1-1, to=2-1]
	\arrow["{(\cdot t^{p-1}) \otimes \id}", from=1-2, to=2-2]
	\arrow["\ad(\theta)"', from=2-1, to=2-2]
\end{tikzcd}
	\end{center}

Let $I^\bullet$ be a resolution of $\ad(E_\star) ^{p,-p}$ and $J^\bullet$ be a resolution of $\widehat{\ad(E_\star)^{p-1,1-p}} \otimes \Omega_X^1(\log D)$. Then,
\begin{align*}
	\mathbb{H}^{p,1-p} = H^1\left(\Gamma\left(I^0\right) \lra \Gamma\left(I^1\right) \oplus \Gamma\left(J^0\right) \lra \Gamma\left(I^2\right) \oplus \Gamma\left(J^1\right) \lra \cdots \right),
\end{align*}
and we can compute the $\bc^\times$-action on $\mathbb{H}^{p,1-p}$ using the above description. Here, we are implicitly using the isomorphism in \Cref{rem: tangent_space_adjoint}.

The complex 
\begin{center}
	\begin{tikzcd}
	{\ad(E_\star)^{p,-p}} & {\widehat{\ad(E_\star)^{p-1,1-p}}\otimes \Omega_X^1(\log D)}
	\arrow["{t\ad(\theta)}", from=1-1, to=1-2]
\end{tikzcd}
\end{center}
is resolved by ${I^\bullet} \to {J^\bullet}'$, where ${J^\bullet}'$ and $J^\bullet$ are related by multiplication by $t$. This is because we have an isomorphism of complexes
\begin{center}
\begin{tikzcd}
	{\ad(E_\star)^{p,-p}} & {\widehat{\ad(E_\star)^{p-1,1-p}}\otimes \Omega_X^1(\log D)} \\
	{\ad(E_\star)^{p,-p}} & {\widehat{\ad(E_\star)^{p-1,1-p}}\otimes \Omega_X^1(\log D)}\rlap{.}
	\arrow["{{\ad(\theta)}}", from=1-1, to=1-2]
	\arrow["\id"', from=1-1, to=2-1]
	\arrow["{{\cdot t}}", from=1-2, to=2-2]
	\arrow["{t\ad(\theta)}"', from=2-1, to=2-2]
\end{tikzcd}
\end{center}
To obtain the chain map, we can compose the vertical maps to get
\begin{center}
\begin{tikzcd}
	{\Gamma\left(I^0\right)} & {\Gamma\left(I^1\right)\oplus \Gamma\left(J^0\right)} & {\Gamma\left(I^2\right) \oplus \Gamma\left(J^1\right)} & \cdots \\
	{\Gamma\left(I^0\right)} & {\Gamma\left(I^1\right)\oplus \Gamma\left({J^0}'\right)} & {\Gamma\left(I^2\right) \oplus \Gamma\left({J^1}'\right)} & \cdots \\
	{\Gamma\left(I^0\right)} & {\Gamma\left(I^1\right)\oplus \Gamma\left(J^0\right)} & {\Gamma\left(I^2\right) \oplus \Gamma\left(J^1\right)} & \cdots\rlap{.}
	\arrow[from=1-1, to=1-2]
	\arrow["{{t^{p}}}"', from=1-1, to=2-1]
	\arrow[from=1-2, to=1-3]
	\arrow["{{(\cdot t^{p})\oplus (\cdot  t^{p-1})}}", from=1-2, to=2-2]
	\arrow[from=1-3, to=1-4]
	\arrow["{(\cdot t^p) \oplus (\cdot t^{p - 1})}", from=1-3, to=2-3]
	\arrow[from=2-1, to=2-2]
	\arrow["\id"', from=2-1, to=3-1]
	\arrow[from=2-2, to=2-3]
	\arrow["{{\id \oplus (\cdot t)}}", from=2-2, to=3-2]
	\arrow[from=2-3, to=2-4]
	\arrow["{\id \oplus (\cdot t)}", from=2-3, to=3-3]
	\arrow[from=3-1, to=3-2]
	\arrow[from=3-2, to=3-3]
	\arrow[from=3-3, to=3-4]
\end{tikzcd}
\end{center}

The chain map is given by multiplication by $t^p$, so $\bc^\times$ acts on $\mathbb{H}^{p,1-p}$ by multiplication by $t^p$.
\end{proof}

\Cref{thm: deform_minimal_energy} says that the minimal energy Higgs bundles always exist when the moduli space $M$ is smooth, and a generic parabolic Higgs bundle in $M$ can be deformed to one of minimal energy. We prove this by using the parabolic Hitchin map defined in \cite[Section~5, p.~495]{Yokogawa_compact_moduli_space} as well as Bialynicki-Birula theory. We briefly review the definitions and properties here.

\subsection{The Hitchin map}

\begin{definition}\label{def: definition_parabolic_hitchin_map}
The \emph{parabolic Hitchin map} is the map $h\colon M \to \oplus_{j = 1}^r H^0(X, \Omega_X^1(\log D)^{\otimes j})$ defined by
\begin{align*}
	h((E_\star, \theta)) = \left(\tr\left( \wedge^1 \theta\right), \tr\left( \wedge^2 \theta\right), \dots, \tr\left(\wedge^r \theta\right)\right).
\end{align*}
The map $h$ sends a parabolic Higgs bundle to the coefficients of the characteristic polynomial of its Higgs field.
\end{definition}

\begin{proposition}\label{prop: parabolic_hitchin_map_proper}
The parabolic Hitchin map is proper.
\end{proposition}
\begin{proof}
	See \cite[Theorem 5.10]{Yokogawa_compact_moduli_space}, where they prove that the valuative criterion for properness holds.
\end{proof}

\begin{lemma}\label{lem: parabolic_deform_C-VHS}
Let $(E_\star, \theta)$ be a stable strongly parabolic Higgs bundle on $X$. Then, the limit $\lim_{t \to 0} (E_\star, t \theta)$ exists $($in the coarse moduli space $M)$ and underlies a $\bc$-VHS. In particular, $(E_\star, \theta)$ can be deformed $($in the coarse moduli space $M)$ to a semi-stable parabolic Higgs bundle underlying a $\bc$-VHS.
\end{lemma}
\begin{proof}
This is proven in \cite[Proposition 1.9(3)]{Mochizuki_Kobayashi_Hitchin}.
\end{proof}

\subsection{The Bialynicki-Birula stratification}

Let $G = \bc^\times$, $Z$ be a some smooth variety over $\bc$ with an action of $G$ on $Z$, and $Z^G$ the locus of fixed points of the $G$-action. Suppose $Z$ can be covered by $G$-invariant quasi-affine opens. Then, Bialynicki-Birula \cite[Theorems 4.1 and~4.3]{Bialynicki_Birula} proves the following. 

\begin{theorem}\label{thm: Bialynicki-Birula}
Let $Z^G = \cup_{i = 1}^m (Z^G)_i$ be the decomposition of $Z^G$ into connected components. Then for any $i = 1,\dots, m$, there are unique locally closed nonsingular $G$-invariant subschemes $Z_i^+$ and $Z_i^{-}$ and unique morphisms $\gamma_i^+ \colon Z_i^+ \to (Z^G)_i$ and $\gamma_i^- \colon Z_i^- \to (Z^G)_i$ such that
\begin{enumerate}[label={\rm(\alph*)}, ref={\rm\alph*}]
	\item\label{t:BB-a} $(Z^G)_i$ is a closed subscheme of\, $Z_i^+$ $($resp.\ $Z_i^-)$ and $\gamma_i^+|_{(Z^G)_i}$ $($resp.\ $\gamma_i^-|_{(Z^G)_i})$ is the identity;
	\item\label{t:BB-b} $Z_i^+$ $($resp.\ $Z_i^-)$ with its induced action of\, $G$ and with the map $\gamma_i^+$ $($resp.\ $\gamma_i^-)$ is a $G$-fibration over $(Z^G)_i$;
	\item\label{t:BB-c} for any closed $a \in (Z^G)_i$, we have the equalities of tangent spaces $T_a(Z_i^+) = T_a(Z)^0 \oplus T_a(Z)^+$ and $T_a(Z_i^-) = T_a(Z)^0 \oplus T_a(Z)^-$, where for a $G$-module $V$, $V^0$ is the $G$-invariant submodule of\, $V$, $V^+$ is spanned by the elements $v \in V$ such that $t \cdot v = t^k v$ for $k > 0$, and $V^-$ is spanned by the elements $v \in V$ such that $t \cdot v = t^k v$ for $k < 0$.
	\end{enumerate}
	If furthermore $\lim_{t \to 0} z \cdot t$ exists for all $z \in Z$ $($resp.\ $\lim_{t \to \infty} z \cdot t$ exists for all $z \in Z)$, then
	\begin{enumerate}[label={\rm(\alph*)}, ref={\rm\alph*},resume]
	  	\item\label{t:BB-d} $\cup Z_i^+ = Z$ $($resp.\ $\cup Z_i^{-} = Z)$;
		\item\label{t:BB-e} $(Z_i^+)^G = (Z^G)_i$ $($resp.\ $(Z_i^-)^G = (Z^G)_i)$ for $i = 1,\dots, m$;
		\item\label{t:BB-f} $Z_i^+$ $($resp.\ $Z_i^-)$ with the map $\gamma_i^+$ $($resp.\ $\gamma_i^-)$ is a $G$-fibration over $(Z^G)_i$; 
		\item\label{t:BB-g} for any closed $a \in (Z^G)_i$, $T_a(Z_i^+) = T_a((Z^G)_i) \oplus T_a(Z_i)^+$ $($resp.\ $T_a(Z_i^-) = T_a((Z^G)_i) \oplus T_a(Z_i)^-)$.
	\end{enumerate}
\end{theorem}

\begin{remark}\label{rem: Completeness_Bialynicki-Birula}
The statements of \Cref{thm: Bialynicki-Birula}(\ref{t:BB-d}), (\ref{t:BB-e}), (\ref{t:BB-f}), and (\ref{t:BB-g}) in \cite[Theorem 4.3]{Bialynicki_Birula} require that our variety $Z$ is complete. This is to ensure that $\lim_{t \to 0} t \cdot z$ exists, so requiring the weaker assumption that $\lim_{t \to 0} t \cdot z$ exists is sufficient. Here, (\ref{t:BB-f}) follows from (\ref{t:BB-b}) and (\ref{t:BB-e}). Similarly, (\ref{t:BB-g}) is the same as (\ref{t:BB-c}) except that it says that the $G$-invariant vectors in $T_a(Z_i)$ are tangent to the $G$-fixed locus $(Z_i)^G$ at $a$.
\end{remark}

\Cref{thm: Bialynicki-Birula} gives us a stratification of $M$ (when it is smooth) compatible with the $\bc^\times$-action on $M$. This is because smoothness of $M$ over $\bc$ ensures that $M$ is a normal variety with a $\bc^\times$-action, so by \cite[Corollary 2]{Sumihiro_Equivariant_Completion}, we can cover $M$ with $\bc^\times$-invariant affine open subsets of $M$. By \Cref{lem: parabolic_deform_C-VHS}, $\lim_{t \to 0} (E_\star, t \theta)$ exists for all $(E_\star, \theta)$, so parts (\ref{t:BB-d}), (\ref{t:BB-e}), (\ref{t:BB-f}), and (\ref{t:BB-g}) of \Cref{thm: Bialynicki-Birula} also apply. We can realize $M$ as a fibration over the locus of $\bc^\times$-fixed points, which by \cite[Theorem 8]{Simpson_harmonic_bundles} is precisely the $\bc$-VHS locus.

\subsection{Proof of \Cref{thm: deform_minimal_energy}}
We are now ready to prove that if $M$ is smooth, then minimal energy parabolic Higgs bundles exist. By \Cref{prop: Non-Abelian_Hodge_Theorem}, this is the same as showing that the minimal energy local systems in $\chi(\pi_1(X \setminus D), C)$ exist.

\begin{proof}[Proof of \Cref{thm: deform_minimal_energy}]

We use the Bialynicki-Birula stratification (\Cref{thm: Bialynicki-Birula}). By assumption, $M$ is smooth, and \Cref{lem: parabolic_deform_C-VHS} implies that $\lim_{t \to 0} (E_\star, t \theta)$ exists, so all parts of \Cref{thm: Bialynicki-Birula} hold in our setting. We first prove the theorem when $M$ is connected.

The $\bc^\times$-fixed locus inside $M$ is precisely the $\bc$-VHS locus by \cite[Theorem 8]{Simpson_harmonic_bundles}. Let $Y$ be this locus and $Y = Y_1 \cup \dots \cup Y_m$ its decomposition into connected components. By \Cref{thm: Bialynicki-Birula}, there are locally closed nonsingular $\bc^\times$-invariant subschemes $M_i^+$ and morphisms $\gamma_i^+\colon M_i^+ \to Y_i$ such that $Y_i$ is inside $M_i^+$ and $\gamma_i^+$ restricts to the identity on $Y_i$. Assume without loss of generality that $M_1^+$ has the greatest dimension among the $M_i^+$ for $i = 1,\dots, m$. Note that $\dim M_1^+ = \dim M$. 

Let $y \in Y_1$. We claim that $y$ is of minimal energy. We define the following subspaces of $T_{y}M$:
\begin{align*}
	W^- &= \spn\left\{ v \in T_{y}M : t \cdot v = t^k v \text{ for some }k < 0\right\},\\
	W^0 &= \spn\left\{ v \in T_{y}M : t \cdot v = v\right\},\\
	W^+ &= \spn\left\{ v \in T_{y}M : t \cdot v = t^k v \text{ for some }k > 0\right\}.
\end{align*}
By \Cref{prop: action_on_tangent_space}, $W^- \cong \oplus_{p < 0} \mathbb{H}^{p,1-p}$, $W^0 \cong \mathbb{H}^{0,1}$, and $W^+ \cong \oplus_{p > 0} \mathbb{H}^{p,1-p}$. Since $\dim_\bc \mathbb{H}^{p,1-p} = \dim_\bc \mathbb{H}^{1-p, p}$ for all $p \in \bz$ (see \Cref{rem: definition_adjoint}), to show that $\mathbb{H}^{p, 1-p} = \mathbb{H}^{1-p, p} = 0$ for $p > 1$, it is enough to show that $W^- = 0$ since 
\begin{align*}
\dim W^- = \dim \left ( \bigoplus_{p < 0} \mathbb{H}^{p,1-p}\right ) = \dim\left (\bigoplus_{p > 1} \mathbb{H}^{p, 1-p}\right ).
\end{align*}
By parts (\ref{t:BB-c}) and (\ref{t:BB-g}) of \Cref{thm: Bialynicki-Birula}, we know that $T_{y}M_1^+ = W^0 \oplus W^+$. But since $\dim M_1^+ = \dim M$, this implies that
\begin{align*}
	W^0 \oplus W^+ = T_{y}M_1^+ = T_{y}M = W^- \oplus W^0 \oplus W^+.
\end{align*}
But then, $W^- = 0$ so $y$ is indeed of minimal energy.

Then to obtain the statement when $M$ is not connected, we observe that the above argument works when restricting to a connected component of $M$.
\end{proof}

\subsection{Some properties of minimal energy local systems}

The minimal energy local systems make up ``most'' of the $\bc$-VHS in $M$. We make this precise below.

\begin{proposition}\label{prop: top_dimension_C-VHS}
Suppose $M$ is smooth, and let $Y$ be the locus of\, $\bc$-VHS inside a connected component of\, $M$. A $\bc$-VHS in $Y$ is of minimal energy if and only if it lies in a connected component of\, $Y$ of largest dimension. The $($complex$)$ dimension of this connected component is half of the $($complex$)$ dimension of\, $M$.
\end{proposition}
\begin{proof}
By \Cref{prop: action_on_tangent_space,rem: tangent_space_adjoint}, the tangent space of $M$ at a $\bc$-VHS $y \in M$ is of the form $T_{y}M = \mathbb{H}^{1}(\ad(E_\star) \to \widehat{ \ad(E_\star)} \otimes \Omega_X^1(\log D)) = \oplus_p \mathbb{H}^{p,1-p}$, and the action of $\bc^\times$ on $\mathbb{H}^{p,1-p}$ is given by multiplication by $t^p$. Hence, parts (\ref{t:BB-c}) and (\ref{t:BB-g}) of \Cref{thm: Bialynicki-Birula} imply that $\mathbb{H}^{0,1}$ is isomorphic $T_{Y,p}$ (which is the tangent space of $Y$ at $p$).

 Let $\widetilde{Y} \subseteq Y$ be the connected component of $Y$ containing $y$. Since $\mathbb{H}^{p,1-p} \cong \mathbb{H}^{1-p,p}$ as vector spaces for all $p \in \bz$, this implies that
\begin{align*}
	\dim \widetilde{Y} = \dim T_{y}Y = \dim \mathbb{H}^{0,1} \leq \frac{1}{2} \dim T_{y}M = \frac{1}{2}\dim M,
\end{align*} 
with equality if and only if $y$ is of minimal energy. Since we know by \Cref{thm: deform_minimal_energy} that minimal energy $\bc$-VHS always exist, we know that the top-dimensional component of $Y$ is always of dimension half of that of~$M$ (and necessarily must be the locus of minimal energy $\bc$-VHS).

If $q = (F_\star, \omega) \in Y$ lies in a component of smaller dimension $Y' \subseteq Y$, then
\begin{align*}
	\dim Y' = \dim T_{q}Y = \dim \mathbb{H}^{0,1}\left(\ad(F_\star) \to \widehat{ \ad(F_\star)} \otimes \Omega_X^1(\log D)\right)  < \frac{1}{2} \dim T_{q}M = \frac{1}{2}\dim M, 
\end{align*}
which implies that there is some $k > 1$ such that $\mathbb{H}^{k, 1-k}(\ad(F_\star) \to \widehat{ \ad(F_\star)}\otimes \Omega_X^1(\log D))$ is nonzero, which violates the definition of minimal energy.
\end{proof}

If $X$ is a positive-genus curve, then the minimal energy parabolic Higgs bundles must correspond to unitary representations of $\pi_1(X \setminus D)$. We make this precise below in \Cref{prop: minimal_energy_genus_0}. 

\begin{proposition}\label{prop: minimal_energy_genus_0}
Let $X$ be a smooth proper curve of genus $g > 0$ and $D$ a reduced effective divisor. If\, $(E_\star, \theta)$ is a stable parabolic Higgs bundle of minimal energy $($with parabolic divisor $D)$ on $X$ of parabolic degree zero, then the Higgs field is zero $($and so $(E_\star, \theta)$ comes from a unitary representation$)$.
\end{proposition}
\begin{proof}
We write $E_\star = E^N_\star \oplus \dots \oplus E^1_\star$, where the graded pieces of our Higgs fields $\theta^p\colon E^p \to E^{p - 1} \otimes \Omega_X^1(\log D)$ are all nonzero. We know by \Cref{prop: Non-Abelian_Hodge_Theorem} that $N = 1$ if and only if $(E_\star, \theta)$ corresponds to a unitary local system. We assume that $N > 1$ for the sake of contradiction. Since $(E_\star, \theta)$ is stable, it is a smooth point of the moduli space of semi-stable parabolic Higgs bundles of parabolic degree zero on $X$ with parabolic divisor~$D$ and parabolic weights $\{\alpha_i^j\}$. Therefore by \Cref{rem: tangent_space_adjoint}, the tangent space to $M$ at $(E_\star, \theta)$ is isomorphic to $\mathbb{H}^1(\ad(E_\star) \to \widehat{\ad(E_\star)} \otimes \Omega_X^1(\log D))$.
	
	If $\mathbb{H}^1(\ad(E_\star) \to \widehat{\ad(E_\star)} \otimes \Omega_X^1(\log D))$ is of Hodge length $1$, we have that
	\begin{align*}
		\mathbb{H}^1\left(0 \to \widehat{\shom(E^1_\star, E^N_\star)_\star} \otimes \Omega_X^1(\log D)\right) = \mathbb{H}^1\left(\shom(E^N, E^1)_\star \to 0\right) = \mathbb{H}^{1-N, N} = 0. 
	\end{align*}
	Then, the Riemann--Roch theorem for curves yields
	\begin{align*}
		0 &\leq  h^0\left(X, \shom\left(E^N_\star, E^1_\star\right)_\star\right) - h^1\left(X, \shom\left(E^N_\star, E^1_\star\right)_\star\right) \\
		&= \deg \shom\left(E^N_\star, E^1_\star\right)_\star + \rank \shom\left(E^N_\star, E^1_\star\right)_\star (1 - g).
	\end{align*}
	Therefore,
	\begin{align*}
		1 - g &\geq  -\frac{ \deg \shom\left(E^N_\star, E^1_\star\right)_\star}{\rank \shom\left(E^N_\star, E^1_\star\right)_\star}.
	\end{align*}
	Then, we have that 
	\begin{align*}
		\pardeg \shom\left(E_\star^N, E^1_\star\right)_\star = \pardeg E^1_\star \cdot \rank E^N_\star - \pardeg E^N_\star \cdot \rank E^1_\star < 0
	\end{align*}
	because $\pardeg E^1_\star < 0$ and $\pardeg E_\star^N > 0$ since $(E_\star, \theta)$ is a stable parabolic Higgs bundle. Then, $\deg \shom(E_\star^N, E^1_\star)_\star \leq \pardeg \shom(E_\star^N, E^1_\star)_\star < 0$, so
	\begin{align*}
		1 - g \geq -\frac{ \deg \shom\left(E^N_\star, E^1_\star\right)_\star}{\rank \shom\left(E^N_\star, E^1_\star\right)_\star} > 0.
	\end{align*}
	Therefore, $1 > g$. But we assumed that $g \geq 1$, so our assumption that $N > 1$ must be false. Therefore, $N = 1$, which implies that $(E_\star, \theta)$ corresponds to a unitary local system under nonabelian Hodge theory (\Cref{prop: Non-Abelian_Hodge_Theorem}).
	\end{proof}

The author thanks Daniel Litt for sharing the following argument learned from a conversation with Bertrand Deroin and Nicolas Tholozan.

\begin{proposition} \label{prop: compact_component}
Let $X$ be a smooth proper curve and $D$ a reduced effective divisor on $X$. Let $(E_\star, \theta) = (\oplus E^k_\star, \oplus \theta_k)$ be a minimal energy parabolic Higgs bundle of rank $r$ on $X$ with parabolic divisor $D$ and generic parabolic weights $($in the sense of \Cref{def: definition_generic_distinct}\,$)$. If\, $p = \sum \rank_{k \in \bz} E^{2k}$ and $q = \sum_{k \in \bz} \rank E^{2k + 1}$, then $p + q = r$ and the connected component of the locus of variations of Hodge structure containing $(E_\star, \theta)$ is compact. In particular, this connected component is a compact connected component of $\chi(\pi_1(X \setminus D),\su(p,q),C)$.
\end{proposition}
\begin{proof}
	Let $\mathbb{V}$ be a local system corresponding to a polarizable $\bc$-VHS. Then, the representation $\rho$ corresponding to $\mathbb{V}$ preserves the polarization, and the polarization is a type $(p,q)$ Hermitian form. Hence, $\rho$ is an $\su(p,q)$ representation of $\pi_1(X\setminus D)$. The $\bc$-VHS correspond to parabolic Higgs bundles with nilpotent Higgs field. Therefore, the parabolic Higgs bundles coming from $\bc$-VHS lie in $h^{-1}(0)$, where $h$ is the parabolic Hitchin map (see \Cref{def: definition_parabolic_hitchin_map}).

	We now consider the minimal energy local systems. As minimal energy local systems underlie a $\bc$-VHS, they necessarily correspond to $\su(p,q)$-representations and also lie in $h^{-1}(0)$. \Cref{prop: parabolic_hitchin_map_proper} tells us that~$h$ is proper, so $h^{-1}(0)$ is compact.
	
	 Let $\chi(\pi_1(\Sigma_{g,n}), \su(p,q), C)$ be the $\su(p,q)$-relative character variety.  The real algebraic group $\su(p,q)$ is a real form of $\sl_r(\bc)$, so $\su(p,q) \times_{\spec(\br)} \spec(\bc) = \sl_r(\bc)$. Therefore, $$\chi(\pi_1(\Sigma_{g,n}), \su(p,q), C) \times_{\spec(\br)} \spec(\bc) = \chi(\pi_1(\Sigma_{g,n}), \sl_r(\bc), C).$$ Hence, $\chi(\pi_1(\Sigma_{g,n}), \su(p,q), C)$ is closed inside the $\gl_r(\bc)$-relative character variety (and it consists of the fixed points of the involution coming from complexifying).
	
	Let $M_{\su(p,q)} \subseteq M$ be the locus of Higgs bundles in $M$ corresponding to $\chi(\pi_1(\Sigma_{g,n}), \su(p,q), C)$ under nonabelian Hodge theory, and let $Y \subseteq M_{\su(p,q)} \subseteq M$ be the locus of minimal energy Higgs bundles. Since $\dim Y = \frac{1}{2} \dim M = \dim M_{\su(p,q)}$, we know that $Y$ is a connected component of $M_{\su(p,q)}$. Therefore, it is closed in $M$. But furthermore, $Y$ consists \emph{entirely} of $\bc$-VHS, implying that it is a subset of $h^{-1}(0)$. Therefore, $Y$ is a component of $M_{\su(p,q)} \cap h^{-1}(0)$. By the properness of $h^{-1}(0)$, we know that $ M_{\su(p,q)} \cap h^{-1}(0)$ is compact since $M_{\su(p,q)}$ is closed in $M$ and closed subsets of compact sets are compact. Finally, $Y$ is a connected component of a compact set $M_{\su(p,q)} \cap h^{-1}(0)$, so $Y$ itself is compact.
	\end{proof}

This mirrors the lower-rank case in \cite[Proposition 2.6]{Deroin_Tholozan_super_maximal}, where the Deroin--Tholozan representations (which are minimal energy local systems in rank $2$) form a compact component of the real character variety.

\begin{remark}\label{rem:old examples}
  The compact components found in Deroin--Tholozan, Tholozan--Toulisse,  and Feng--Zhang are examples of minimal energy components. This is because all of the local systems found in those components underlie $\bc$-VHS, which by \Cref{prop: top_dimension_C-VHS} is equivalent to being of minimal energy.
\end{remark}

\section{Genus zero}\label{sec:genus zero}

We restrict our attention to minimal energy parabolic Higgs bundles on $\bp^1$. In contrast to the higher-genus case (see \Cref{prop: minimal_energy_genus_0}), where the minimal energy local systems correspond to unitary representations, the minimal energy local systems do not necessarily underlie a $\bc$-VHS with a trivial Hodge filtration. However, we show in \Cref{thm: bounds_C-VHS} that when the number of punctures on $\bp^1$ is large compared to the rank of the minimal energy local system, the minimal energy energy local systems underlie a $\bc$-VHS with at most two steps in their Hodge filtration.

Our ultimate goal in this section is to prove \Cref{thm: bounds_C-VHS}. To that end, let $(E_\star, \theta)$ be a parabolic Higgs bundle on $\bp^1$ with distinct parabolic weights and parabolic divisor $D$ corresponding to an irreducible minimal energy local system. Since $(E_\star, \theta)$ underlies a $\bc$-VHS, it is of the form $(\oplus_{p = 1}^N E^p_\star, \oplus \theta^p)$, where our Higgs field $\theta^p \colon E^p_\star \to E^{p - 1}_\star \otimes \Omega_{\bp^1}^1(\log D)$ is nonzero (due to the irreducibility assumption on the corresponding local system) and $N$ is the number of graded pieces. The number of steps in the Hodge filtration of the $\bc$-VHS corresponding to $(E_\star, \theta)$ is the same as the number of graded pieces of $(E_\star, \theta)$. 

As an intermediate step, we will show that if $N > 2$, then the number of punctures $n$ cannot be too large compared to $\rank E_\star$. Let $V^k_\star = \oplus_{p} \shom(E_\star^p, E_\star^{p + k})_\star = \ad(E_\star)^{k, -k}$, and let $\ad(\theta)^{k} \colon V^k_\star \to V^{k - 1}_\star \otimes \Omega_{\bp^1}^1(\log D)$ be the induced Higgs field on $\ad(E_\star)$.

\begin{proposition}\label{prop: main_bound}
	Let $(E_\star, \theta)$ be a minimal energy parabolic Higgs bundle on $(\bp^1,x_1 + \dots + x_n)$ with $N$ graded pieces and with generic unitary local monodromy in the sense of \Cref{def: definition_generic_distinct} $($so that there are $r$ distinct parabolic weights of\, $(E_\star, \theta))$.  If\, $N \geq 3$, then
\begin{align*}
	n \leq \frac{2 \rank V^{N-1} - \rank V^1 + 2 \sum_{k = 2}^{N-1}\rank V^{N-1} + \sum_{k = 2}^{N-1}(2k - 3) \rank V^k}{\sum_{k = 2}^{N - 1}\left( \rank V^{k - 1} - \rank \coker \ad(\theta)^{k}\right)  + \sum_{k = 2}^{N - 1}(k - 2) \rank V^k}.
\end{align*}
\end{proposition}

\begin{remark}
	The denominator in this expression can be simplified, but we will leave it unsimplified because in \Cref{lem: bounds_$1$} we will explicitly bound the term $\sum_{k = 2}^{N - 1}( \rank V^{k - 1} - \rank \coker \ad(\theta)^{k})$ by a positive number. In particular, we will show the denominator is always nonzero. We also note that the term $\sum_{k = 2}^{N - 1}(k - 2) \rank V^k$ only contributes when $N > 3$, and when $N = 3$, the sum vanishes. We do not index starting from $k = 3$ to avoid specifying what an empty sum means. 
\end{remark}

The bound in \Cref{prop: main_bound} depends on $(E_\star, \theta)$ itself and not just the rank of the parabolic Higgs bundle. To extract a bound from \Cref{prop: main_bound} depending only on $\rank E_\star$, we bound each of the individual terms $V^k$ in terms of $\rank E_\star$. As a consequence, whenever the rank of $E_\star$ violates the given bound, we know that $N \leq 2$, proving \Cref{thm: bounds_C-VHS}.

To that end, we prove a series of lemmas that place restrictions on the behavior of the Higgs field of $(\ad(E_\star), \oplus \ad(\theta)^{k})$ whenever $(E_\star, \theta)$ is of minimal energy.

\begin{lemma}\label{lem: ker_coker_vanish}
If\, $(E_\star, \theta)$ is a parabolic Higgs bundle of minimal energy with distinct parabolic weights, then for all $1 < k < N$, 
\begin{align*}
  H^1\left(\bp^1, \ker \ad(\theta)^{k}\right) = H^0\left(\bp^1, \coker \ad(\theta)^{k}\right) = 0.
\end{align*}
\end{lemma}

\begin{proof}
	By \Cref{def: definition_minimal_energy} and \Cref{rem: tangent_space_adjoint}, $\mathbb{H}^1(\ad(E_\star) \to \widehat{\ad(E_\star)} \otimes \Omega_{\bp^1}^1 (\log D))$ is of Hodge length $1$. Note that since our parabolic weights are distinct, $\widehat{V^k_\star} \otimes \Omega_{\bp^1}^1(\log D) = V^k_\star \otimes \Omega_{\bp^1}^1(\log D)$ for all $1 \leq k \leq N$. Therefore, the complex $\ad(E_\star) \to \widehat{\ad(E_\star)} \otimes \Omega_{\bp^1}^1 (\log D)$ has associated graded pieces
	\begin{align*}
		0 & \lra V^{N - 1}_\star \otimes \Omega_{\bp^1}^1(\log D)\\
		V^{N-1}_\star & \lra V^{N - 2}_\star \otimes \Omega_{\bp^1}^1(\log D)\\
		&\hspace{0.75em}\vdots\\
		 V^{2 - N}_\star & \lra V^{1 - N}_\star \otimes \Omega_{\bp^1}^1(\log D)\\
		V^{1 - N}_\star & \lra 0.
	\end{align*}
	The graded pieces of $\mathbb{H}^1(\ad(E_\star) \to \widehat{\ad(E_\star)} \otimes \Omega_{\bp^1}^1 (\log D))$ are given by
	\begin{align*}
		\mathbb{H}^{k,1-k}\left(\ad(E_\star) \to \widehat{\ad(E_\star)} \otimes \Omega_{\bp^1}^1(\log D)\right) = \mathbb{H}^1\left(V^{k}_\star \to V^{k - 1}_\star \otimes \Omega_{\bp^1}^1 (\log D)\right)
	\end{align*}
	whenever $k \neq 1$. For $k = 1$,
	\begin{align*}
		\mathbb{H}^{1,0}\left(\ad(E_\star) \to \widehat{\ad(E_\star)} \otimes \Omega_{\bp^1}^1(\log D)\right) = \mathbb{H}^1\left(V^{1}_\star \to \widehat{V^{0}_\star} \otimes \Omega_{\bp^1}^1 (\log D)\right)
	\end{align*}
	but $\widehat{V^0_\star} \otimes \Omega_{\bp^1}^1(\log D) \neq V^0 \otimes \Omega_{\bp^1}^1(\log D)$. We have an exact sequence of complexes of sheaves
	\begin{align*}
		0 \lra \left(\ker \ad(\theta)^{k}\right)[0] \lra \left[V^{k}_\star \to V^{k - 1}_\star \otimes \Omega_{\bp^1}^1 (\log D)\right] \lra \left(\coker \ad(\theta)^{k}\right)[-1] \lra 0, 
	\end{align*}
	so we get a long exact sequence on hypercohomology	
	\begin{center}
	\begin{tikzcd}
	0 \\
	{\mathbb{H}^0\left(\bp^1, \left(\ker \ad(\theta)^{k}\right)[0]\right)} & {\mathbb{H}^0\left(V^{k}_\star \to V^{k - 1}_\star \otimes \Omega_{\bp^1}^1 (\log D)\right)}
        & {\mathbb{H}^0\left(\bp^1, \left(\coker \ad(\theta)^{k}\right)[-1]\right)} \\
	{\mathbb{H}^1\left(\bp^1, \left(\ker \ad(\theta)^{k}\right)[0]\right)} & {\mathbb{H}^1\left(V^{k}_\star \to V^{k - 1}_\star \otimes \Omega_{\bp^1}^1 (\log D)\right)} & {\mathbb{H}^1\left(\bp^1, \left(\coker \ad(\theta)^{k}\right)[-1]]\right)} \\
	&& 0\rlap{.}
	\arrow[from=1-1, to=2-1]
	\arrow[from=2-1, to=2-2]
	\arrow[from=2-2, to=2-3]
	\arrow[from=2-3, to=3-1]
	\arrow[from=3-1, to=3-2]
	\arrow[from=3-2, to=3-3]
	\arrow[from=3-3, to=4-3]
\end{tikzcd}
	\end{center}
	Since
        \begin{align*}
        \mathbb{H}^0\left(\bp^1, \coker \ad(\theta)^{k}[-1]\right) &= 0,\\
        \mathbb{H}^1\left(\bp^1, \coker \ad(\theta)^{k}[-1]\right) &= H^0\left(\bp^1, \coker \ad(\theta)^{k}\right),\\
        \intertext{and} 
        \mathbb{H}^1\left(\bp^1, \ker \ad(\theta)^{k}[0]\right) &= H^1\left(\bp^1, \ker \ad(\theta)^{k}\right),
        \end{align*}
        we have a short exact sequence
	\begin{align*}
		0 \lra H^1\left(\bp^1, \ker \ad(\theta)^{k}\right) \lra \mathbb{H}^1\left(V^{k}_\star \to V^{k - 1}_\star \otimes \Omega_{\bp^1}^1 (\log D)\right) \lra H^0\left(\bp^1, \coker \ad(\theta)^{k}\right) \lra 0.
	\end{align*}
	By definition, $\mathbb{H}^{k, 1 - k}(\ad (E_\star) \to \widehat{\ad (E_\star)}\otimes \Omega_{\bp^1}^1(\log D)) = \mathbb{H}^1(V^k_\star \to V^{k - 1}_\star\otimes \Omega_{\bp^1}^1(\log D))$. But since $(E_\star, \theta)$ is of minimal energy, $\mathbb{H}^{k, 1 - k}(\ad (E_\star) \to \widehat{\ad (E_\star)}\otimes \Omega_{\bp^1}^1(\log D)) = 0$ for $1 < k < N$. Therefore, the terms $H^1(\bp^1, \ker \ad(\theta)^{k})$ and $H^0(\bp^1, \coker \ad(\theta)^{k})$ vanish for $1 < k < N$ due to the above sequence being short exact and the vanishing of the term in the middle.
\end{proof}

Lemma~\ref{lem: ker_coker_vanish} gives us strong conditions on the Higgs fields of minimal energy Higgs bundles.

\begin{corollary}\label{cor: coker_locallyfree}
Let $(E_\star, \theta)$ be a minimal energy parabolic Higgs bundle on $\bp^1$. Then for $1 < k < N$, $\coker \ad(\theta)^{k}$ is a locally free sheaf.
\end{corollary}
\begin{proof}
	As $\bp^1$ is a smooth curve and $\coker \ad(\theta)^{k}$ is a coherent sheaf, we obtain a splitting
	$\coker \ad(\theta)^{k} = (\coker \ad(\theta)^{k} )_{\lf} \oplus (\coker \ad(\theta)^{k} )_{\tors}$, where $(\coker \ad(\theta)^{k} )_{\lf}$ is locally free and $(\coker \ad(\theta)^{k} )_{\tors}$ is a torsion sheaf. Since $(E_\star, \theta)$ is a minimal energy parabolic Higgs bundle, $\mathbb{H}^{1}(\ad(E_\star) \to \widehat{\ad(E_\star)} \otimes \Omega_{\bp^1}^1(\log D))$ is of Hodge length $1$. \Cref{lem: ker_coker_vanish} implies that for $1 < k < N$,
	\begin{align*}
		0 = H^0\left(\bp^1, \coker \ad(\theta)^{k}\right) = H^0\left(\bp^1, \left(\coker \ad(\theta)^{k}\right)_{\lf}\right) \oplus H^0\left(\bp^1, \left(\coker \ad(\theta)^{k} \right)_{\tors}\right).
	\end{align*}
	But $(\coker \ad(\theta)^{k} )_{\tors}$ is a torsion sheaf on a curve and is hence supported at finitely many points; see \cite[Exercise II.5.6(c)]{Hartshorne}. Therefore, $H^0(\bp^1, (\coker \ad(\theta)^{k} )_{\tors}) = 0$ if and only if $(\coker \ad(\theta)^{k} )_{\tors}$ has empty support, which holds if and only if $(\coker \ad(\theta)^{k} )_{\tors} = 0$. Therefore, $\coker \ad(\theta)^{k} = (\coker \ad(\theta)^{k} )_{\lf}$, as desired.
\end{proof}

\begin{corollary}\label{cor: bounds_ker_coker}
If $(E_\star, \theta)$ is a minimal energy parabolic Higgs bundle on $\bp^1$, then for $1 < k < N$, $\deg \coker \ad(\theta)^{k} \leq - \rank \coker \ad(\theta)^{k}$ and $\deg \ker \ad(\theta)^{k} \geq -\rank \ker \ad(\theta)^{k}$.
\end{corollary}
\begin{proof}
	As $(E_\star, \theta)$ is a minimal energy Higgs bundle on $\bp^1$, by \Cref{lem: ker_coker_vanish}, whenever $1 < k < N$, we have that $H^0(\bp^1, \coker \ad(\theta)^{k}) = H^1(\bp^1, \ker \ad(\theta)^{k}) = 0$.
	
	As $\ker \ad(\theta)^{k}$ is a vector bundle on $\bp^1$, $\ker \ad(\theta)^{k} \cong \mathscr{O}(a_1) \oplus \dots \oplus \mathscr{O}(a_{\rank \ker \ad(\theta)^{k}})$, where the $a_j$ are integers and $a_1 + \dots + a_{\rank \ker \ad(\theta)^{k}} = \deg \ker \ad(\theta)^{k}$. Therefore, 
	\begin{align*}
		H^1\left(\bp^1, \ker \ad(\theta)^{k}\right) = H^1\left(\bp^1, \mathscr{O}\left(a_1\right)\right) \oplus \dots \oplus H^1\left(\bp^1, \mathscr{O}\left(a_{\rank \ker \ad(\theta)^{k}}\right)\right) = 0
	\end{align*}
	if and only if all of the integers $a_j$ are at least $-1$. So, $\deg \ker \ad(\theta)^{k} = a_1 + \dots + a_{\rank \ker \ad(\theta)^{k}} \geq - \rank \ker \ad(\theta)^{k}$.
	
	By  \Cref{cor: coker_locallyfree}, $\coker \ad(\theta)^{k}$ is a vector bundle and is isomorphic to $\mathscr{O}(b_1) \oplus \dots \oplus \mathscr{O}(b_{\rank \coker \ad(\theta)^{k}})$. Since
	\begin{align*}
		H^0\left(\bp^1, \coker \ad(\theta)^{k}\right) = H^0\left(\bp^1, \mathscr{O}\left(b_1\right)\right) \oplus \dots \oplus H^0\left(\bp^1, \mathscr{O}\left(b_{\rank \coker \ad(\theta)^{k}}\right)\right) = 0,
	\end{align*}
	we know that all of the integers $b_j$ are at most $-1$. Therefore,
\begin{equation*}\pushQED{\qed}
		\deg \coker \ad(\theta)^{k} = b_1 + \dots + b_{\rank \coker \ad(\theta)^{k}} \leq -\rank \coker \ad(\theta)^{k}.\qedhere \popQED
	\end{equation*}
\renewcommand{\qed}{} 
\end{proof}

We now prove a series of technical lemmas used in the proof of \Cref{prop: main_bound}.

\begin{lemma} \label{lem: bounds_$1$}
If\, $(E_\star, \theta)$ is of minimal energy on $\bp^1$ with distinct parabolic weights, then 
\begin{align*}
\pardeg V^1)\star & \geq (\deg D - 1) \sum_{k = 2}^{N - 1}(k - 2)\rank V^k - \sum_{k = 2}^{N -1} (k - 1)\rank V^k
\end{align*}
and $\deg V^{N - 1} = \rank V^{N - 1}\cdot(1 - \deg D)$.
\end{lemma}
\begin{proof}
We first show that $\deg V^{N - 1} = \rank V^{N- 1} \cdot (1 - \deg D)$. Since $(E_\star, \theta)$ is of minimal energy, $\mathbb{H}^1(V^{1 - N} \to 0) = H^1(\bp^1, V^{1 - N}) = 0$. But observe that $V^{1 -N }_\star$ is a sub-parabolic Higgs bundle of $\ad (E_\star)$ with zero Higgs field, so $V^{1 - N}$ and all of its subbundles are also sub-Higgs bundles. Therefore, all subbundles of $V^{1-N}$ must have nonpositive parabolic degree. Since all of the parabolic weights on $V^{1 - N}$ are nonzero (due to the hypothesis that the parabolic weights are distinct), this implies that every subbundle of $V^{1-N}$ must have \emph{negative} degree. Therefore, as vector bundles, $V^{1 -N} \cong \mathscr{O}(a_1) \oplus \dots \oplus \mathscr{O}(a_{\rank V^{1-N}})$ with $a_j < 0$ for all $j = 1,\dots, \rank V^{1-N}$. But since $H^1(\bp^1, V^{1 - N}) = 0$, we have that $a_j \geq - 1$ for all $j = 1,\dots, \rank V^{1-N}$. Therefore, $V^{1 - N} \cong \mathscr{O}(-1)^{\oplus \rank V^{1 - N}}$. Equivalently, $V^{N - 1} \cong (\mathscr{O}(1)^{\oplus\rank V^{N - 1}})(- D)$ and so $\deg V^{N - 1} = \rank V^{N - 1}\cdot (1 - \deg D)$.
	
	\vspace{.5cm}
	
	We now prove the other part of the claim and show that
	\begin{align*}
\pardeg V^1_\star & \geq (\deg D - 1) \sum_{k = 2}^{N - 1}(k - 2)\rank V^k - \sum_{k = 2}^{N -1} (k - 1)\rank V^k.
\end{align*}

	We have exact sequences
	\begin{align*}
		0 \lra \ker \left(\ad(\theta)^{k}\right)_\star \lra V^k_\star \lra V^{k - 1}_\star \otimes \Omega_{\bp^1}^1(\log D) \lra \coker \left(\ad(\theta)^{k}\right)_\star \lra 0.
	\end{align*}
	 Since degrees are additive in exact sequences,
	 \begin{align*}
	 	\deg V^{k - 1} \otimes \Omega_{\bp^1}^1(\log D) + \deg \ker \ad(\theta)^{k}  = \deg V^k + \deg \coker \ad(\theta)^{k} 
	 \end{align*}
	 for all $2 < k < N$. Then using $\deg V^{k - 1} \otimes \Omega_{\bp^1}^1(\log D) = \deg V^{k - 1} + \rank V^{k - 1}\cdot (\deg D - 2)$, we have that  
	 \begin{align*}
	 	\deg V^{k - 1} &= \deg V^k + \rank V^{k - 1} \cdot (2 - \deg D)  + \deg \coker \ad(\theta)^{k}  - \deg \ker \ad(\theta)^{k}.
	 \end{align*}
	 
	We have that $\deg \coker \ad(\theta)^{k} \leq - \rank \coker \ad(\theta)^{k}$ and $\deg \ker \ad(\theta)^{k} \geq - \rank \ker \ad(\theta)^{k}$ by \Cref{cor: bounds_ker_coker}. Therefore, 
	\begin{align*}
		\deg V^{k - 1} &\leq  \deg V^k + \rank V^{k - 1} \cdot (2 - \deg D) + \rank \ker \ad(\theta)^{k} - \rank \coker \ad(\theta)^{k}.
	\end{align*}
	Since $\rank \ker \ad(\theta)^{k} - \rank \coker \ad(\theta)^{k} = \rank V^k - \rank V^{k - 1}$, we then have that
	\begin{align*}
		\deg V^{k - 1} & \leq \deg V^k + \rank V^{k - 1} \cdot (1 - \deg D) + \rank V^{k}.
	\end{align*}
		By induction, for $\ell > 0$ we obtain the bound
	\begin{align*}
		\deg V^{k - \ell} & \leq \deg V^k + (1 - \deg D) \cdot  \sum_{j = 1}^\ell \rank V^{k - j} + \sum_{j = 1}^{\ell} \rank V^{k - j + 1}.
	\end{align*}
	In the special case where $k = N - 1$,  by substituting $\deg V^{N - 1} = \rank V^{N - 1}\cdot (1 - \deg D)$, we have  that 
	\begin{align*}
		\deg V^{N - 1 - \ell} & \leq (1 - \deg D) \cdot  \sum_{j = 0}^\ell \rank V^{N - 1 - j} + \sum_{j = 1}^{\ell} \rank V^{N - j}.
	\end{align*}

	We now observe that $V^0_\star \oplus V^{-1}_\star \oplus \dots \oplus V^{1 - N}_\star$ is a sub-Higgs bundle of $\ad(E_\star)$, so by the semi-stability of $\ad (E_\star)_\star$, we have that $\pardeg (V^0_\star\oplus V^{-1}_\star \oplus \dots \oplus V^{1 - N}_\star) \leq 0$. Equivalently, $\pardeg(V^{N - 1}_\star \oplus \dots \oplus V^1_\star) \geq 0$. We use this to first bound $\pardeg V^{k}_\star$ when $k \geq 2$.
		
	By applying \Cref{lem: bounds_pardeg} to $\pardeg V_\star^k$, we see that
	\begin{align*}
		\pardeg V^k_\star \leq \deg V^k + \rank V^k \cdot \deg D.
	\end{align*}
	Then substituting our bounds  $\deg V^k$ obtained above, we have that 
	\begin{align*}
		\pardeg V^k_\star & \leq \sum_{j = k}^{N - 1} \rank V^j + (1 - \deg D) \cdot \sum_{j = k + 1}^{N - 1} \rank V^j.
	\end{align*}

Since $\pardeg(V_\star^{N - 1} \oplus \dots \oplus V_\star^1) \geq 0$, we have that 
\begin{align*}
	\pardeg V^{1}_\star &\geq - \sum_{k = 2}^{N - 1} \pardeg V^k_\star\\
	& \geq - \sum_{k = 2}^{N - 1} \left (\sum_{j = k}^{N - 1} \rank V^j + (1 - \deg D) \cdot \sum_{j = k + 1}^{N - 1} \rank V^j \right ),
\end{align*}
so we have obtained a lower bound in terms of the $\rank V^k$. Simplifying yields
\begin{equation*}\pushQED{\qed}
	\pardeg V^1_\star \geq  (\deg D - 1) \sum_{k = 2}^{N - 1}(k - 2)\rank V^k - \sum_{k = 2}^{N -1} (k - 1)\rank V^k.\qedhere \popQED
	\end{equation*}
\renewcommand{\qed}{}     
\end{proof}

\begin{lemma}\label{lem: bounds_rank_coker}
For all $1 < k < N$, we have that  
\begin{align*}
	\rank V^{k - 1} - \rank \coker \ad(\theta)^{k} \geq N - k.
\end{align*}
In particular, $\rank V^{k - 1} - \rank \coker \ad(\theta)^{k} > 0$.
\end{lemma}
\begin{proof}
	Since $V^{k - 1} \otimes \Omega_{\bp^1}^1(\log D)$ surjects onto $\coker \ad(\theta)^{k}$, we know that $\rank V^{k - 1} - \rank \coker \ad(\theta)^{k} \geq 0$. By \Cref{cor: bounds_ker_coker}, we have that 
	\begin{align*}
		\rank V^{k - 1} - \rank \coker \ad(\theta)^{k} = \rank V^k - \rank \ker \ad(\theta)^{k},
	\end{align*}
	and hence if $\rank V^{k - 1} - \rank \coker \ad(\theta)^{k} = 0$, then $\ad(\theta)^{k}$ must be the zero map. We show that $\rank V^k - \rank \ker \ad(\theta)^{k} \geq N - k$. Since $\ker \ad(\theta)^{k}$ is a vector bundle, we can compute its rank by computing the rank of a fiber. Equivalently, we can compute the dimension of $(\ad(\theta)^{k})_s$ for any point $s \in \bp^1$. We do this by explicitly writing out the map $\ad(\theta)^{k}$.

	The map 
	\begin{align*}
		\ad(\theta)^{k}\colon \oplus_{j = 1}^{N - k}\shom\left(E^j_\star, E^{j + k}_\star\right)_\star \lra \oplus_{j = 1}^{N - k + 1}\shom\left(E^j_\star, E^{j + k-1}_\star\right)_\star \otimes \Omega_{\bp^1}^1(\log D)
	\end{align*}
	is given by
	\begin{align*}
		\ad(\theta)^{k}\left(f_1^{1 + k}, \dots, f_{N - k}^N\right)
		 &=  \left( \theta^{1 + k} \circ f_1^{1 + k},  \theta^{k + 2} \circ f_2^{k + 2} -\left(f_1^{1 + k} \otimes \id\right) \circ \theta^2, \dots,\right. \\
		 & \left.\qquad \theta^N \circ f_{N - k}^N - \left(f_{N - k - 1}^{N - 1}\otimes \id\right)\circ \theta^{N -k}, - \left(f_{N - k}^k \otimes \id \right) \circ \theta^{N - k + 1}\right),
	\end{align*}
	where $f_j^{j + k} \in \shom(E^j_\star, E^{j + k}_\star)_\star$.
	
	Let $s \in \bp^1$ be any point for which none of the Higgs fields of $\theta^p \colon E^p \to E^{p - 1} \otimes \Omega_{\bp^1}(\log D)$ vanishes. Such a point exists because our Higgs fields are nonzero, and the locus where each $\theta^p$ vanishes is a closed subvariety of $\bp^1$ of positive codimension. By \Cref{cor: coker_locallyfree}, we know that $\coker \ad(\theta)^{k}$ is a locally free sheaf for all $1 < k < N$. Therefore, $\rank \ker \ad(\theta)^{k} = \dim \ker (\ad(\theta)^{k})_s$. Then, $\ad(\theta)^{k}$ at $s$ is given by
	\begin{align*}
		(\ad(\theta)^{k})_s &= \begin{pmatrix}
			 B_{2}^{k + 1} & A_1^{k} &  &   & \\
			  &  B_3^{k + 2} & A_2^{k + 1} & \\
			  & & \ddots & \ddots  \\
			   & &   &  B_{N - k + 1}^r &  A_{N - k}^{N - 1}  \\
		\end{pmatrix},
	\end{align*}
	where 
	\begin{align*}
	A_j^{j + k - 1}	\left(f_j^{j + k}\right) &= \left(\theta^{j + k}\right)_s \circ f_j^{j + k}\\
	B_{j + 1}^{j + k}\left(f_{j}^{j + k}\right) &= -\left(f_{j}^{j + k} \otimes \id\right) \circ \left(\theta^{j + 1}\right)_s.
	\end{align*}
	Then,
	\begin{align*}
		 \rank \left(\ad(\theta)^{k}\right)_s \geq \sum_{j = 1}^{N - k} \rank \left(A_j^{j + k - 1}\right),
	\end{align*}
	and since $s$ is a general point, we know that our Higgs fields $\theta^p \colon E^{p} \to E^{p - 1} \otimes \Omega_{\bp^1}^1(\log D)$ are nonzero, and so our maps $A_j^{j +k - 1}$ are also nonzero. Therefore, since $\rank (A_j^{j + k - 1}) \geq 1$, we have that
	\begin{align*}\pushQED{\qed}
		\rank V^{k - 1} - \rank \coker \ad(\theta)^{k} &= \rank V^k - \rank \ker\ad(\theta)^{k} = \dim \left(V^k\right)_s - \dim \ker \left(\ad(\theta)^{k}\right)_s \\
		& = \rank \left(\ad(\theta)^{k}\right)_s \\
		&\geq N - k.\qedhere \popQED
	\end{align*}
\renewcommand{\qed}{}   
\end{proof}

We can now prove \Cref{prop: main_bound}.

\begin{proof}[Proof of \Cref{prop: main_bound}]

We recall what we want to prove. Let $(E_\star, \theta)$ be a minimal energy parabolic Higgs bundle on $(\bp^1,D)$ with $N$ graded pieces and with generic unitary local monodromy (in the sense of \Cref{def: definition_generic_distinct}). Let $D = x_1 + \dots + x_n$. If $N \geq 3$, then
\begin{align*}
	\deg D = n \leq \frac{2 \rank V^{N-1} - \rank V^1 + 2 \sum_{k = 2}^{N-1}\rank V^{N-1} + \sum_{k = 2}^{N-1}(2k - 3) \rank V^k}{\sum_{k = 2}^{N - 1}\left( \rank V^{k - 1} - \rank \coker \ad(\theta)^{k}\right)  + \sum_{k = 2}^{N - 1}(k - 2) \rank V^k}.
\end{align*}

We have an exact sequence of parabolic bundles
\begin{align*}
	0 \lra \ker \left(\ad(\theta)^{k}\right)_\star \lra V^k_\star \lra V^{k - 1}_\star \otimes \Omega_{\bp^1}^1(\log D) \lra \coker \left(\ad(\theta)^{k}\right)_\star \lra 0,
\end{align*}
and by  \Cref{cor: coker_locallyfree}, for $1 < k < N$, this is an exact sequence of vector bundles on $\bp^1$. Therefore,  by \Cref{lem: par_deg_properties}, we obtain that 
\begin{align*}
 \pardeg \ker \left(\ad(\theta)^{k}\right)_\star + \pardeg V^{k - 1}_\star \otimes \Omega_{\bp^1}^1(\log D) = \pardeg V^k_\star + \pardeg \coker \left(\ad(\theta)^{k}\right)_\star.
\end{align*}
Since $\pardeg V^{k - 1}_\star \otimes \Omega_{\bp^1}^1(\log D) = \pardeg V^{k - 1}_\star + \rank V^{k - 1} \cdot ( \deg D - 2)$, by summing the above equation over $1 < k < N$, we obtain the equality
\begin{align*}
		&\sum_{k = 2}^{N - 1}\pardeg \ker\left(\ad(\theta)^{k}\right)_\star + \sum_{k = 2}^{N - 1}\pardeg V^{k - 1}_\star + (\deg D - 2) \sum_{k = 2}^{N - 1}\rank V^{k - 1}\\
		&= \sum_{k = 2}^{N - 1} \pardeg V^{k}_\star + \sum_{k = 2}^{N - 1} \pardeg \coker\left(\ad(\theta)^{k}\right)_\star.
\end{align*}
We have a telescoping sum that simplifies to
\begin{align*}
	&\sum_{k = 2}^{N - 1}\pardeg \ker\left(\ad(\theta)^{k}\right)_\star + \pardeg V^{1}_\star + (\deg D - 2) \sum_{k = 2}^{N - 1}\rank V^{k - 1}\\
		&= \pardeg V^{N-1}_\star + \sum_{k = 2}^{N - 1} \pardeg \coker\left(\ad(\theta)^{k}\right)_\star.
\end{align*}
The combination of \Cref{lem: bounds_$1$,lem: bounds_pardeg} yields the inequalities 
\begin{align*}
		\pardeg V^1_\star &\geq  (\deg D - 1) \sum_{k = 2}^{N - 1}(k - 2)\rank V^k - \sum_{k = 2}^{N -1} (k - 1)\rank V^k,\\
		\pardeg V^{N-1}_\star &\leq \deg V^{N-1} + \deg D \cdot \rank V^{N - 1}  = \rank V^{N - 1}.
\end{align*}
Then we use the above inequalities. After simplifying and rearranging our original expression, we obtain the inequality
\begin{align*}
	&\sum_{k = 2}^{N - 1} \pardeg \coker\left(\ad(\theta)^{k}\right)_\star - \sum_{k = 2}^{N - 1}\pardeg \ker\left(\ad(\theta)^{k}\right)_\star \\
	& \geq (\deg D - 2) \sum_{k = 2}^{N - 1} \rank V^{k - 1} + (\deg D - 1) \sum_{k = 2}^{N - 1}(k - 2)\rank V^k - \sum_{k = 2}^{N -1} (k - 1)\rank V^k - \rank V^{N - 1}.
\end{align*}
We now apply \Cref{lem: bounds_pardeg} to the parabolic bundles $\ker (\ad(\theta)^{k})_\star$ and $\coker (\ad(\theta)^{k})_\star$. By \Cref{lem: bounds_pardeg}, we have that for all $1 < k < N$, 
\begin{align*}
	\pardeg \coker \left(\ad(\theta)^{k}\right)_\star &\leq \deg \coker \left(\ad(\theta)^{k}\right) + \rank \coker \left(\ad(\theta)^{k}\right) \cdot \deg D,\\
		\deg \ker\left(\ad(\theta)^{k}\right) &\leq \pardeg \ker \left(\ad(\theta)^{k}\right)_\star.
\end{align*}
Then, combining the inequalities yields 
\begin{align*}
 &(\deg D - 2) \sum_{k = 2}^{N - 1} \rank V^{k - 1} + (\deg D - 1) \sum_{k = 2}^{N - 1}(k - 2)\rank V^k - \sum_{k = 2}^{N -1} (k - 1)\rank V^k - \rank V^{N - 1}\\	& \leq \sum_{k = 2}^{N - 1} \deg \coker(\ad(\theta)^{k}) - \sum_{k = 2}^{N - 1}\deg \ker(\ad(\theta)^{k}) + \deg D \cdot \left (\sum_{k = 2}^{N - 1}\rank \coker( \ad(\theta)^{k}) \right ).
\end{align*}
Since
\begin{align*}
	\rank \coker \left(\ad(\theta)^{k}\right)- \rank \ker \left(\ad(\theta)^{k}\right) & \leq \deg \ker \left(\ad(\theta)^{k}\right) - \deg \coker \left(\ad(\theta)^{k}\right)
\end{align*}
by \Cref{cor: bounds_ker_coker}, we have that
\begin{align*}
  & (\deg D - 2) \sum_{k = 2}^{N - 1} \rank V^{k - 1} + (\deg D - 1) \sum_{k = 2}^{N - 1}(k - 2)\rank V^k - \sum_{k = 2}^{N -1} (k - 1)\rank V^k - \rank V^{N - 1} \\
  & \leq \sum_{k = 2}^{N - 1} \left(\rank \ker \left(\ad(\theta)^{k}\right) - \rank \coker \left(\ad(\theta)^{k}\right)\right) + \deg D \cdot \left (\sum_{k = 2}^{N - 1}\rank \coker \left(\ad(\theta)^{k}\right) \right ).
\end{align*}
Then, the rank-nullity theorem implies that $\rank \ker (\ad(\theta)^{k}) - \rank \coker (\ad(\theta)^{k})  = \rank V^k -\rank V^{k - 1}$. So, substituting yields the inequality
\begin{align*}
	&  (\deg D - 2) \sum_{k = 2}^{N - 1} \rank V^{k - 1} + (\deg D - 1) \sum_{k = 2}^{N - 1}(k - 2)\rank V^k - \sum_{k = 2}^{N -1} (k - 1)\rank V^k - \rank V^{N - 1}\\& \leq \sum_{k = 2}^{N - 1} (\rank V^k -\rank V^{k - 1}) + \deg D \cdot \left (\sum_{k = 2}^{N - 1}\rank \coker( \ad(\theta)^{k}) \right ).
\end{align*}
But the expression $\sum_{k = 2}^{N - 1} (\rank V^k -\rank V^{k - 1})$ is a telescoping sum and equals $\rank V^{N-1} - \rank V^1$. Therefore, we obtain the inequality
\begin{align*}
	 &  (\deg D - 2) \sum_{k = 2}^{N - 1} \rank V^{k - 1} + (\deg D - 1) \sum_{k = 2}^{N - 1}(k - 2)\rank V^k - \sum_{k = 2}^{N -1} (k - 1)\rank V^k - \rank V^{N - 1}\\& \leq \rank V^{N-1} - \rank V^1 + \deg D \cdot \left (\sum_{k = 2}^{N - 1}\rank \coker \ad(\theta)^{k}\right ).
\end{align*}
Simplifying yields
\begin{align*}
	&\deg D \cdot \left ( \sum_{k = 2}^{N - 1}\left( \rank V^{k - 1} - \rank \coker \left(\ad(\theta)^{k}\right)\right) + \sum_{k = 2}^{N - 1}(k - 2) \rank V^k\right )\\
	& \leq 	2 \rank V^{N-1} - \rank V^1 + 2 \sum_{k = 2}^{N-1}\rank V^{k-1} + \sum_{k = 2}^{N-1}(2k - 3) \rank V^k.
\end{align*}
Then, \Cref{lem: bounds_rank_coker} tells us that $\rank V^{k - 1} - \coker \ad(\theta)^{k} \geq N - k$ for all $1 < k < N$ and $N \geq 3$. Therefore, the left-hand side is positive, since one of the differences $\rank V^{k - 1} - \coker \ad(\theta)^{k}$ is bounded by a positive integer. Dividing then gives us
\begin{equation*}\pushQED{\qed}
	\deg D \leq \frac{2 \rank V^{N-1} - \rank V^1 + 2 \sum_{k = 2}^{N-1}\rank V^{k-1} + \sum_{k = 2}^{N-1}(2k - 3) \rank V^k}{\sum_{k = 2}^{N - 1}\left( \rank V^{k - 1} - \rank \coker \left(\ad(\theta)^{k}\right)\right)  + \sum_{k = 2}^{N - 1}(k - 2) \rank V^k}. \qedhere \popQED
\end{equation*}
\renewcommand{\qed}{}    
\end{proof}

We can now extract a uniform bound on $\deg D$ in terms of $\rank E_\star$ from \Cref{prop: main_bound}.

\begin{proof}[Proof of \Cref{thm: bounds_C-VHS}]
	If $(E_\star, \theta)$ has more than two nonzero pieces in its grading, $(E_\star, \theta)$ is given by a Higgs bundle $(E_\star, \theta_\star)$ with nonzero Higgs field as it underlies an irreducible nonunitary local system. As our local monodromy data is generic in the sense of \Cref{def: definition_generic_distinct}, our parabolic weights are all distinct. And the $\theta^p$ are never identically zero as $\mathbb{V}$ underlies an irreducible nonunitary representation. Let $N$ be the number of graded pieces in $(E_\star, \theta)$ (\textit{i.e.}, the number of steps in the Hodge filtration of the $\bc$-VHS corresponding to $(E_\star, \theta)$). Suppose that $N > 2$ and $n = \deg D \geq 5r^2 - 25$.
	
	We can apply \Cref{prop: main_bound}. We have that 
	\begin{align*}
	\deg D \leq \frac{2 \rank V^{N-1} - \rank V^1 + 2 \sum_{k = 2}^{N-1}\rank V^{k-1} + \sum_{k = 2}^{N-1}(2k - 3) \rank V^k}{\sum_{k = 2}^{N - 1}\left( \rank V^{k - 1} - \rank \coker \ad(\theta)^{k}\right)  + \sum_{k = 2}^{N - 1}(k - 2) \rank V^k}.
	\end{align*}
	
By \Cref{lem: bounds_rank_coker}, we have that $\rank V^{k - 1} - \rank \coker \ad(\theta)^{k} \geq N - k$. Therefore, 
\begin{align*}
	\deg D \leq \frac{2 \rank V^{N-1} - \rank V^1 + 2 \sum_{k = 2}^{N-1}\rank V^{k-1} + \sum_{k = 2}^{N-1}(2k - 3) \rank V^k}{\sum_{k = 2}^{N - 1}(N - k)  + \sum_{k = 2}^{N - 1}(k - 2) \rank V^k}.
\end{align*}
	The denominator is bounded below by
	\begin{align*}
		\sum_{k = 2}^{N - 1}(N - k) = \frac{(N - 1)(N - 2)}{2},
	\end{align*}
	and this is nonzero (and positive) whenever $N > 2$. Then, we coarsely bound the numerator from above. Since $\rank V^j \geq 1$, we know that 
	\begin{align*}
	\rank \left(\oplus_{j = 1 - N, j \neq k}^{N - 1} V^j\right) = \sum_{j = 1 - N, j \neq k}^{N - 1} \rank V^j \geq \sum_{j = 1 - N, j \neq k}^{N - 1} 1 = 2(N-1).
	\end{align*}

Since $\rank \ad(E_\star) = r^2 - 1$, we have that $\rank V^k = \rank \ad (E_\star) - \rank (\oplus_{j = 1 - N, j \neq k}^{N - 1} V^j) \leq  r^2 - 1 - 2(N - 1) = r^2 - 2N + 1$. We now bound the numerator. We have that 
	
\begin{align*}
&2 \rank V^{N-1} - \rank V^1 + 2 \sum_{k = 2}^{N-1}\rank V^{k-1} + \sum_{k = 2}^{N-1}(2k - 3) \rank V^k\\
	&\leq  2  (r^2 - 2N + 1) - 1 + 2 \sum_{k = 2}^{N-1}(r^2 - 2N + 1) + \sum_{k = 2}^{N-1}(2k - 3)(r^2 - 2N + 1)\\
	&\leq (r^2 - 2N + 1)(N^2 - 2N + 2) - 1\\
	& \leq (r^2 - 2N + 1)(N^2 - 2N + 2)\\
	& \leq (r^2 - 5)(N^2 - 2N + 2).
	\end{align*}
	Since the denominator is at least $\frac{(N - 1)(N - 2)}{2}$,  by \Cref{prop: main_bound}, we have that
	\begin{align*}
		\deg D \leq \frac{2(r^2 - 5)(N^2 - 2N + 2)}{(N-1)(N-2)}.
	\end{align*}
	Then for $N \geq 3$, we have that
	\begin{align*}
		\frac{N^2 - 2N + 2}{(N-1)(N - 2)} \leq \frac{5}{2}.
	\end{align*}
	Therefore,
	\begin{align*}
		\deg D\leq \frac{2(r^2 - 5)(N^2 - 2N + 2)}{(N-1)(N-2)} \leq 5(r^2 - 5) = 5r^2 - 25.
	\end{align*}
	But we assumed that $\deg D > 5r^2 - 25$, so our assumption that $N > 2$ was false. Therefore, $N \leq 2$.
	\end{proof}

\section{Consequences of Theorem~\ref{thm: bounds_C-VHS}}\label{sec:consequences}

\Cref{thm: bounds_C-VHS} states that when $\deg D > 5r^2 - 25$,  a stable minimal energy parabolic Higgs bundle of rank~$n$ has at most two graded pieces. We show that this is the best we can do. We construct an example of a stable minimal energy parabolic Higgs bundle on $\bp^1$ in every rank with parabolic divisor of arbitrarily large degree that has two graded pieces. We first prove a technical lemma.

\begin{lemma}\label{lem: technical_lemma_for_example}
Let $a$ and $r$ be positive integers with $r \geq 2$, and let $0 < \epsilon < \frac{1}{1 + ra}\left ( \frac{r-1}{r}\right )$ be a very small real number $($depending on $r$ and $a)$. Then,
\begin{align*}
	\frac{1}{r-1}\left (1 - \frac{1 +a}{1 + ra} + \epsilon  \right ) < \frac{1 + a}{1 + ra} - \epsilon.
\end{align*}
\end{lemma}
\begin{proof}
We have that
\begin{align*}
	\epsilon < \frac{1}{1 + ra}\left ( \frac{r-1}{r}\right ).
\end{align*}
Then,
\begin{align*}
	\frac{1}{1 + ra}\left ( \frac{r-1}{r}\right ) = \frac{r + ra - 1 - ra}{r(1 + ra)} = \frac{r(1 + a) - (1 + ra)}{r(1 + ra)} = \frac{1 + a}{1 + ra} - \frac{1}{r}.
\end{align*}
Therefore,
\begin{align*}
	\epsilon < \frac{1 + a}{1 + ra} - \frac{1}{r}.
\end{align*}
So,
\begin{align*}
	1 + r \epsilon  < r \cdot \frac{1 + a}{1 + ra} .
\end{align*}
Then adding $-\frac{1 + a}{1 + ra}  - (r - 1)\epsilon$ to both sides yields 
\begin{align*}
1 - \frac{1 + a}{1 + ra}  + \epsilon< 	(r - 1) \left (\frac{1 + a}{1 + ra}  - \epsilon \right ),
\end{align*}
giving the desired inequality
\begin{equation*}\pushQED{\qed}
 \frac{1}{r - 1}\left ( 1 - \frac{1 + a}{1 + ra}  + \epsilon \right ) < 	 \frac{1 + a}{1 + ra}  - \epsilon .\qedhere \popQED
	\end{equation*}
\renewcommand{\qed}{} 
\end{proof}

In the following example, we give an example of a minimal energy local system on $\bp^1 \setminus \{x_1,\dots, x_n\}$ in every rank (with the number of punctures depending on the rank) that has $\su(1,r-1)$-monodromy. These representations $\rho$ admit a $\rho$-equivariant holomorphic map from the universal cover of $\bp^1 \setminus \{x_1,\dots, x_n\}$ to complex hyperbolic space (which is the symmetric space of $\su(1,r-1)$) and have been previously constructed by Tholozan--Toulisse \cite[Theorem 2]{Tholozan_Toulisse}.

\begin{example}\label{ex: two_graded_pieces_example}
Let $a$ be some nonnegative integer, $r > 2$, and $\epsilon < \frac{1}{1 + ra}\left ( \frac{r - 1}{r} \right )$ very small. Let $E = S \oplus Q$, where $S = \mathscr{O}(-a)^{\oplus (r-1)}$ and $Q = \mathscr{O}(-a - 1)$, so $\rank E = r$. We let $D$ be a reduced effective divisor on $\bp^1$ such that $\deg D = 1 + ra$ and $D = x_1 + \dots + x_{1 + ra}$ are the points of the divisor. We choose our Higgs field $\theta \colon S \to Q \otimes \Omega_X^1(\log D)$ to be generic, and since $Q$ is a line bundle, we know that $\theta$ is surjective. At each $x_j$, we give $E$ the weights  
\begin{align*}
	0 < \frac{1}{r-1}\left (1 - \frac{1 +a}{1 + ra} + \epsilon  \right ) + \epsilon_1  < \dots < \frac{1}{r-1}\left (1 - \frac{1 +a}{1 + ra} + \epsilon  \right ) + \epsilon_{r - 1} <\frac{1 + a}{1 + ra} - \epsilon < 1,
\end{align*}
where $\epsilon_1 < \dots < \epsilon_{r -1}$ are chosen to be very close to zero and such that $\epsilon_1 + \dots + \epsilon_{r - 1} = 0$. By \Cref{lem: technical_lemma_for_example}, we can pick $\epsilon_1, \dots, \epsilon_{r - 1}$ satisfying these conditions.

We give $E_{x_j}$ the parabolic flag $ E_{x_j}=S_j^1 \oplus Q_{x_j}\supset \dots \supset S_j^{r - 1}\oplus Q_{x_j}\supset Q_{x_j} \supseteq 0$, where $\{S_j^i\}$ is a generically chosen flag of $S_{x_j}$. Our parabolic weights are distinct, and our parabolic structure can be chosen to be generic in the sense of \Cref{def: definition_generic_distinct} by modifying the real numbers $\epsilon$ and $\epsilon_1,\dots, \epsilon_{r - 1}$. 

We have a nonzero Higgs field because $\shom(S_\star, Q_\star)_\star \otimes \Omega_X^1(\log D) \cong \mathscr{O}(ra)^{\oplus (r-1)}$, which has global sections since $r > 2$ and $a \geq 0$.

Since $(E_\star, \theta)$ only has two graded pieces, the complex $\ad(E_\star) \to \widehat{\ad(E_\star)}\otimes \Omega_X^1(\log D)$ has graded pieces
\begin{align*}
	0 &\lra \shom(Q_\star, S_\star)_\star \otimes \Omega_X^1(\log D),\\
	\shom(Q_\star, S_\star)_\star &\lra \widehat{\shom(S_\star, S_\star)_\star} \otimes \Omega_X^1(\log D) \oplus \widehat{\shom(Q_\star, Q_\star)_\star}\otimes \Omega_X^1(\log D),\\
	\shom(S_\star, S_\star)_\star \oplus \shom(Q_\star, Q_\star)_\star &\lra \shom(S_\star, Q_\star)_\star \otimes \Omega_X^1(\log D),\\
	\shom(S_\star, Q_\star)_\star & \lra 0.
\end{align*}
By the symmetry of the Hodge structure on $\mathbb{H}^1(\ad(E_\star) \to \widehat{\ad(E_\star)}\otimes \Omega_X^1(\log D))$, to show that $E_\star$ is of minimal energy, we need to show that $\mathbb{H}^1(\shom(S_\star, Q_\star)_\star \to 0) = H^1(\bp^1, \shom(S_\star, Q_\star)_\star) = 0$. But this is true, since $\shom(S_\star, Q_\star)_\star \cong \mathscr{O}(-1)^{\oplus(r-1)}$. 

Therefore, it remains to show that $(E_\star,\theta)$ is a stable parabolic Higgs bundle. By construction, the parabolic degree of $E_\star$ is zero. The sub-Higgs bundles of $(E_\star, \theta)$ are the bundles $Q_\star$, $\ker(\theta)_\star$, and $V_\star \oplus Q_\star$, where $V \subseteq S$ is any subbundle.

We first note that $\mu_\star(Q_\star) = \deg Q + \sum_{x_j \in D} (\frac{1 + a}{1 + ra} - \epsilon) = -1 - a + 1 + a - (1 + ra) \epsilon < 0$ since $\deg D = 1 + ra$. Therefore, $Q_\star$ is not destabilizing. 

To check that $\ker(\theta)_\star$ is not destabilizing, we note that it fits in a exact sequence
\begin{align*}
	0 \lra \ker(\theta)_\star \lra S_\star \lra Q_\star \otimes \Omega_X^1(\log D) \lra 0.
\end{align*}
Hence by \Cref{lem: par_deg_properties}, we have the equality
\begin{align*}
	\pardeg \ker (\theta)_\star = \pardeg S_\star - \pardeg Q_\star \otimes \Omega_X^1(\log D).
\end{align*}
Observe that $\pardeg S_\star = \epsilon\deg D$ and $\pardeg Q_\star = -\epsilon\deg D$, so
\begin{align*}
	\pardeg \ker (\theta)_\star = \pardeg S_\star - \pardeg Q_\star \otimes \Omega_X^1(\log D) = 2 \epsilon\deg D - \deg D + 2.
\end{align*}
Since $\epsilon > 0$ is chosen to be very small, as long as $\deg D \geq 3$ (which holds, for example, as long as $a > 0$ and $r > 0$), we know that $\pardeg \ker (\theta)_\star < 0$, and so $\ker(\theta)_\star$ cannot be destabilizing.

It remains to check that $V_\star \oplus Q_\star$ is not destabilizing for any subbundle $V \subseteq S$. We have the equalities $\pardeg (V_\star \oplus Q_\star) = \pardeg V_\star + \pardeg Q_\star = \pardeg V_\star - \epsilon \deg D$. Let $N = \rank V$. We will compute the parabolic degree of $V_\star$ (which is a rank $N$ parabolic bundle), and we will show that the resulting rank $N + 1$ parabolic bundle $V_\star \oplus Q_\star$ has negative parabolic degree.

Since $V$ is a subbundle of $\mathscr{O}(-a)^{\oplus(r - 1)}$ of rank $N$, we have that $\deg V \leq -aN$. Furthermore, the sum of the parabolic weights of $V$ at every point of the divisor is at most the sum of the $N$ largest parabolic weights at every point $x_i$ in $D$. Therefore,
\begin{align*}
	\pardeg V_\star &\leq -aN + \deg D \cdot \left (  \sum_{j = r - N}^{r-1} \left ( \frac{1}{r - 1} \left (1 - \frac{1 + a}{1 + ra} + \epsilon \right ) + \epsilon_{j} \right ) \right )\\
	&= -aN + \frac{N \deg D}{r - 1} \left ( 1 - \frac{1 + a}{1 + ra} + \epsilon \right ) + \deg D(\epsilon_{r - 1} + \dots + \epsilon_{n-1-N}).
\end{align*}
Simplifying and substituting $\deg D = 1 + ra$ into  the bound yields
\begin{align*}
	\pardeg V_\star &\leq  \frac{N \deg D \epsilon}{r - 1} + \deg D( \epsilon_{r - 1} + \dots + \epsilon_{r - 1 - N}).
\end{align*}
Since $\pardeg Q_\star = - \epsilon\deg D$,
\begin{align*}
	\pardeg (V_\star \oplus Q_\star) \leq \frac{N \deg D \epsilon}{r - 1} + \deg D( \epsilon_{r - 1} + \dots + \epsilon_{r - 1 - N}) - \epsilon \deg D.
\end{align*}
Since $0 \leq N < r - 1$ and $\epsilon_{r - 1} + \dots + \epsilon_{r - 1 - N}$ is chosen to be very small, we have $\pardeg (V_\star \oplus Q_\star) < 0$.

Therefore, every sub-Higgs bundle of $(E_\star, \theta)$ fails to be destabilizing, so $(E_\star, \theta)$ is indeed a stable minimal energy parabolic Higgs bundle of rank $r$ with parabolic divisor $D$ that has two graded pieces.

By choosing $a \gg 0$, we can force $\deg D =1 + ra$ to be larger than $5r^2 - 25$. Therefore, we cannot do better than two graded pieces in \Cref{thm: bounds_C-VHS}. In this example, we picked our parabolic weights so that at each point in $D$, the weights sum to an integer. Therefore, the local monodromy data  not only is unitary but actually lies in $\su(n)$. 
\end{example}

\subsection{An application to Gromov--Witten theory} We give an application of minimal energy local systems to Gromov--Witten invariants.

Given a parabolic bundle $V_\star$ of rank $n$ on $\bp^1$ with parabolic divisor $D$, we can obtain a modified complete parabolic bundle $\widetilde{V}_\star$ with the properties given below. This process is described in \cite[below Theorem~5.1]{Agnihotri_Woodward} and in \cite[Appendix]{Belkale_Unitary_Deligne_Simpson}. 

\begin{proposition}\label{prop: properties_modified_bundle}
Let $V_\star$ be a parabolic bundle and $\widetilde{V}_\star$ its modified parabolic bundle.
\begin{enumerate}
	\item\label{p:pmb-1} The underlying bundle of\, $\widetilde{V_\star}$ is a trivial bundle such that $\rank \widetilde{V} = \rank V$. 
	\item\label{p:pmb-2} If\, $V_\star$ has $0 \leq \alpha_j^1 < \dots < \alpha_j^{n_j} < 1$ at $x_j \in D$, then the weights $\{\widetilde{\alpha_j^i}\}$ of\, $\widetilde{V}_\star$ at $x_j \in D$ no longer necessarily lie within $[0,1)$ but satisfy $\widetilde{\alpha_j^1} \leq \dots \leq \widetilde{\alpha_j^{n}} < \widetilde{\alpha_j^1} + 1$. The weights can be repeated.
	\item\label{p:pmb-3} The parabolic degree of a modified parabolic bundle is defined to be the sum of the degree of the underlying bundle and the weights above every point. Then, $\pardeg (V_\star) = \pardeg(\widetilde{V}_\star)$.
	\item\label{p:pmb-4} There is a bijection between subbundles $W_\star$ of\, $V_\star$ and modified parabolic subbundles $\widetilde{W}_\star$ of\, $\widetilde{V}_\star$. This bijection preserves the parabolic degree. If\, $W_\star$ and $\widetilde{W}_\star$ are sent to each other under this bijection, then $\pardeg W_\star = \pardeg \widetilde{W}_\star$. 
	\item\label{p:pmb-5} As a consequence of \eqref{p:pmb-1}, \eqref{p:pmb-3}, and \eqref{p:pmb-4}, $V_\star$ is parabolic stable if and only if\, $\widetilde{V}_\star$ is parabolic stable.
\end{enumerate}
\end{proposition}

A conjugacy class $C \subseteq \su(r)$ is determined by the eigenvalues of a diagonal representative. If the eigenvalues are given by $(e^{2\pi i \theta_1},\dots, e^{2\pi i \theta_r})$, then we obtain real numbers $\theta_1,\dots, \theta_r$, where each $\theta_i$ is determined up to an integer. Note that  $\theta_1,\dots, \theta_r$ must sum to an integer because matrices in $\su(r)$ have determinant~$1$. We choose  $\theta_1, \dots, \theta_r$ such that $\theta_1 + \dots + \theta_r = 0$ and that (after ordering from largest to smallest) $\theta_1 \geq \dots \geq \theta_r \geq \theta_1 - 1$. This uniquely determines the conjugacy class $C$. Let $\lambda(C) = (\theta_1 ,\dots, \theta_r)$, and for some subset $I \subseteq \{1,\dots, r\}$ of size $k$,  let $\lambda_I(C) = \sum_{j \in I} \theta_j$. 

A subset $I \subseteq \{1,\dots, r\}$ of size $k$ determines a Schubert variety $\Omega_I$ in the Grassmannian $\gr(k,n)$, and we let $\sigma_I$ be the cohomology class associated to $\Omega_I$. When $I_1,\dots, I_n$ are all subsets of $\{1,\dots, r\}$ of equal size, we let $\langle \sigma_{I_1},\dots, \sigma_{I_n} \rangle_d$ be the number of degree $d$ rational curves in $\gr(k,n)$ passing through $\Omega_{I_1}, \dots, \Omega_{I_n}$. The number $\langle \sigma_{I_1},\dots, \sigma_{I_n} \rangle_d$ is a \emph{Gromov--Witten invariant} for $\gr(k,n)$. For details on Schubert varieties and Gromov--Witten invariants, see \cite[Definition 4]{Belkale_Unitary_Deligne_Simpson}.

Then, Biswas \cite{Biswas_parabolic} (when $r = 2$), Belkale \cite[Theorem 7]{Belkale_Unitary_Deligne_Simpson}, and Agnihotri--Woodward \cite[Theorem~3.1]{Agnihotri_Woodward} independently prove the following. 

 \begin{theorem}\label{thm: Belkale_inequality}
	There is an $\su(r)$-local system on $\bp^1 \setminus \{x_1,\dots, x_n\}$ with prescribed $\su(r)$-local monodromy $C_1,\dots, C_n$ if and only if for any $k \in \{1,\dots, r - 1\}$ and any choice of subsets $I_1,\dots, I_n \subseteq \{1,\dots, r\}$ of cardinality $k$, the inequality
	\begin{align*}
		\sum_{j = 1}^n \lambda_{I_j}(C_j) \leq d
	\end{align*}
	holds whenever $\langle \sigma_{I_1},\dots, \sigma_{I_n} \rangle_d = 1$.
 \end{theorem}
 
The above statement is taken from \cite[Theorem 7]{Belkale_Unitary_Deligne_Simpson}.

\begin{corollary}\label{cor: Corollary_Gromov-Witten}
Let $X = \bp^1$, and let $D = x_1 + \dots + x_n$ be a reduced effective divisor. Let $C = (C_1,\dots, C_n)$ be a collection of generic and connected local monodromy data $($in the sense of \Cref{def: definition_generic_distinct}\,$)$. Let $(E_\star, \theta)$ be a minimal energy stable parabolic Higgs bundle that only has two graded pieces and has parabolic weights corresponding to the data $C_1,\dots, C_n$. Then, $E_\star = H_\star \oplus (V/H)_\star$ with Higgs field $\theta \colon H \to (V/H)\otimes \Omega_X^1(\log D)$.

Let $\widetilde{E}_\star$ be the modified parabolic bundle of\, $E_\star$ and $\widetilde{H}_\star$ the modified parabolic bundle of\, $H_\star$. For all $j =1,\dots, n$, we let $I_j \subseteq \{1,\dots, r\}$ be a subset of size $\rank H$ such that for all $i \in I_j$, $\theta_j^i$ is a parabolic weight of\, $\widetilde{H}_\star$. Then, the Gromov--Witten invariant $\langle \sigma_{I_1},\dots, \sigma_{I_n} \rangle_{- \deg \widetilde{H}}$ is nonzero.
\end{corollary}
\begin{proof}
Since $(E_\star, \theta)$ has two graded pieces, it is nonunitary. For a connected $n$-tuple of conjugacy classes $C = (C_1,\dots, C_n)$, the relative character variety $\chi(\pi_1(X \setminus D), \gl_r(\bc), C)$ is connected by \cite[Theorem~5.1.1]{Hausel_Letellier_Rodriguez-Villegas}. Therefore, the locus of minimal energy local systems is unique and does not consist of unitary local systems since $(E_\star, \theta)$ is nonunitary.

Since $E_\star = H_\star \oplus (V/H)_\star$ with $H_\star$ a destabilizing subbundle of $V_\star$, we know that the associated Gromov--Witten invariant $\langle \sigma_{I_1},\dots, \sigma_{I_n}\rangle_{- \deg \widetilde{H}}$ is nonzero. This is because $\widetilde{H}$ is a subbundle of a trivial bundle and hence determines a rational degree $-\deg \widetilde{H}$ curve in $\gr(k,r)$ intersecting the Schubert varieties $\Omega_{I_1}, \dots, \Omega_{I_n}$.
\end{proof} 

\begin{remark}\label{rem: main_theorem_application}
When $\deg D > 5r^2 - 25$, where $r = \rank E$, \Cref{thm: bounds_C-VHS} applies and the minimal energy local systems have the desired number of graded pieces. 
\end{remark}

We give an example of  \Cref{cor: Corollary_Gromov-Witten}.

\begin{example}\label{ex: Example_Gromov-Witten}
Let $(E_\star,\theta)$ be a parabolic Higgs bundle on $\bp^1$ with parabolic divisor $D = x_1 + x_2 + x_3$. We let the underlying bundle be $E =\mathscr{O}(-2) \oplus \mathscr{O}(-1) \oplus \mathscr{O}(-1)$. We set $S = \mathscr{O}(-2)$ and $Q = \mathscr{O}(-1) \oplus \mathscr{O}(-1)$. We let $Q_{x_j} \supseteq Q_j \supseteq 0$ be a generic flag of the fiber $Q_{x_j}$ at each point $x_j$. For $j = 1,2$, we let our parabolic flag at $x_j$ be given by
\begin{align*}
	E_{x_j} = Q_{x_j} \oplus S_{x_j} \supseteq Q_j \oplus S_{x_j} \supseteq S_{x_j} \supseteq 0.
\end{align*}
At $x_3$, we let our parabolic flag be given by
\begin{align*}
	E_{x_3} = Q_{x_3} \oplus S_{x_3} \supseteq Q_{x_3} \supseteq Q_3 \supseteq 0.
\end{align*}
Let $\frac{1}{24} >\epsilon > \frac{1}{48}$ be a real number. We choose our weights at $x_1$ to be $\{\frac{1}{2}, \frac{1}{2} + \epsilon, 1 - \epsilon\}$. The weights at $x_2$ are $\{ \frac{1}{8} - \epsilon, \frac{1}{8}, \frac{3}{4} + \epsilon \}$. The weights at $x_3$ are given by $\{ \frac{1}{3} - \epsilon, \frac{1}{3}, \frac{1}{3} + \epsilon \}$. That is, at $x_1$ and $x_2$, we give $S$ the largest parabolic weights. At $x_3$, we give $S$ the smallest parabolic weights.

We note that $\mu(S_\star) = \frac{1}{12} - \epsilon$. Every subbundle of $Q_\star$ has parabolic degree at most $- \frac{1}{24} - \epsilon$. Therefore, $(E_\star, \theta)$ is parabolic stable.
Then, $\shom(S_\star, Q_\star)_\star \cong \mathscr{O}(-1) \oplus \mathscr{O}(-1)$ and $\shom_\star(S_\star,Q_\star)_\star \otimes \Omega_X^1(\log D) \cong \mathscr{O} \oplus \mathscr{O}$. Therefore, $(E_\star, \theta)$ is a minimal energy parabolic Higgs bundle.

Let $I_1$, $I_2$, and $I_3$ be  subsets of $\{1,2,3\}$ of size $1$ so that $\lambda_{I_j}(C_j)$ gives the weight of $S_\star$ over $x_j$. So, our Schubert varieties are subvarieties of $\bp^2$. Here, $I_1 = I_2 = \{1\}$ and $I_3 = \{3\}$. The classes $\sigma_{I_1}$ and $\sigma_{I_2}$ are the class of a point. The class $\sigma_{I_3}$ is the class of $\bp^2$, and $\deg \widetilde{S} = -1$. Therefore, we get that $\langle \sigma_{I_1}, \sigma_{I_2}, \sigma_{I_3}\rangle_1 \neq 0$. Geometrically, this says that there is at least one degree $1$ rational curve in $\bp^2$ that intersects $\bp^2$ and two generic points in $\bp^2$. Since a degree $1$ rational curve is just a line, this says that there is at least one line passing through two points in $\bp^2$ in generic position. \Cref{cor: Corollary_Gromov-Witten} does not, however, tell us that there is a \textit{unique} line passing through two points in $\bp^2$ in generic position.
\end{example}


\providecommand{\bysame}{\leavevmode\hbox to3em{\hrulefill}\thinspace}

\end{document}